\documentclass[11pt]{amsart}

\usepackage{standard}
\usepackage{cancel}
\usepackage{mathtools}
\usepackage{xcolor}

\newcommand{\realtau}{\mathcal{T}}
\renewcommand{\tau}{\sigma}

\newcommand{\gmetric}{\mathfrak{g}}
\newcommand{\hmetric}{\mathfrak{h}}
\newcommand{\A}{\mathcal{A}}
\newcommand{\floor}[1]{\lfloor #1 \rfloor}
\newcommand{\fracpart}[1]{\{ #1 \}}
\newcommand{\N}{\mathbb{N}}

\newcommand{\boxg}{\Box_{\gmetric}}
\newcommand{\la}{\langle}
\newcommand{\ra}{\rangle}

\newcommand{\dd}{\hbox{ }d}
\DeclareMathOperator{\Diffb}{Diff_{b}}
\DeclareMathOperator{\Diff}{Diff}

\newcommand{\LEonorm}[1]{ \|\partial #1\|_{LE} + \|\la r \ra^{-1}#1\|_{LE}}
\newcommand{\LEnorm}[1]{\sup_{m}\|\la r \ra^{-\frac{1}{2}}#1\|_{L^2(\RR_+ \times A_m)}}

\newcommand{\LEsnorm}[1]{ \sum_m \|\la r \ra^{\frac{1}{2}}#1\|_{L^2(\RR_+ \times A_m)}}
\newcommand{\LE}{{\mathcal {LE}}}
\newcommand{\LEs}{{\mathcal {LE}^*}}
\newcommand{\lesn}[1]{_{\mathcal{LE}^{*,#1}}}
\newcommand{\letn}[1]{_{\mathcal{LE}_\sigma^{#1}}}
\newcommand{\lettn}[1]{_{\mathcal{LE}_{\sigma ,1}^{#1}}}
\newcommand{\len}[1]{_{\mathcal{LE}^{#1}}}
\newcommand{\lrsupo}{_{L_r^2L_\omega^\infty(A_m)}}

\newcommand{\ltwoam}{_{L^2(A_m)}}

\newcommand{\les}{_{\mathcal{LE}^*}}
\newcommand{\vijk}{T^i\Omega^jS^k v}
\newcommand{\gijk}{T^i\Omega^jS^k g}
\newcommand{\ijk}{T^i\Omega^jS^k}

\newcommand{\R}{\mathbb{R}}
\newcommand{\F}{\mathcal{F}}
\newcommand{\M}{\mathcal{M}}
\newcommand{\V}{\mathcal{V}}
\newcommand{\hol}{\mathcal{H}}
\newcommand{\RC}{\mathsf{RC}}
\newcommand{\Rnbar}{\overline{\RR^n}}
\newcommand{\loc}{\text{loc}}

\newcommand{\scc}{\text{sc}}
\newcommand{\res}{\text{res}}

\newcommand{\Hb}{H_b}
\newcommand{\tP}{\widetilde{P}}

\newcommand{\angvar}{\theta}

\newcommand{\loss}{\varpi}
\newcommand{\poly}{\mathbb{C}[\sigma]}

\renewcommand{\Im}[1]{\operatorname{Im}#1}

\DeclareMathOperator{\Id}{Id}

\begin{document}
\title[Generalized Price's law]{Generalized Price's law on fractional-order
  asymptotically flat stationary spacetimes}
\author{Katrina Morgan}
\author{Jared Wunsch}
\address{Department of Mathematics, Northwestern University}
\email{katrina.morgan@northwestern.edu, jwunsch@math.northwestern.edu}

\begin{abstract}
We obtain estimates on the rate of decay of a solution to the wave
equation on a stationary spacetime that tends to Minkowski space at a
rate $O(\lvert x \rvert^{-\kappa}),$ $\kappa \in (1,\infty) \backslash
\mathbb{N}.$  Given suitably smooth and decaying initial data, we show a wave
locally enjoys the decay rate $O(t^{-\kappa-2+\epsilon}).$
\end{abstract}

\maketitle

\section{Introduction}
The goal of this work is to study the relationship among pointwise
decay rates of waves, low- and high-frequency resolvent estimates, and
the large scale behavior of the background geometry. We study
solutions to the wave equation on asymptotically flat stationary
4-dimensional Lorentzian spacetimes with signature (3,1).  The flat
Minkowski metric, which we denote $\mathfrak{m}$, is given by
$\mathfrak{m} = dt^2 - dx_1^2-dx_2^2-dx_3^2$. The spacetimes
considered here are of the form
$\mathfrak{g} = \mathfrak{m} + \mathfrak{h}$ where $\mathfrak{h}$ has
metric coefficients which decay like $r^{-\kappa}$ for some
$\kappa \in (1,\infty)\setminus \mathbb{N}$. We say such a metric
tends toward flat at a rate of $r^{-\kappa}$. We find that, given
sufficiently differentiable and decaying Cauchy data, waves decay
locally at a rate of $t^{-\kappa-2+\epsilon}$.    The main new input to this
decay estimate is a certain resolvent estimate valid uniformly near
zero frequency.   In previous work \cite{Mo:20}, the
first author studied the case when $\kappa \in \N$ and established
$t^{-\kappa-2}$ decay rates when the background geometry
  exhibits spherical symmetry, and $t^{-\kappa-2+\epsilon}$ decay in
  the absence of spherical symmetry.
	
The study of pointwise decay rates on asymptotically flat spacetimes
arises in general relativity. In \cite{Price}, physicist Richard Price
gave a heuristic argument anticipating a $t^{-3}$ pointwise decay rate
for waves on the Schwarz\-schild spacetime, which describes space in
the presence of a single, non-rotating black hole. 
This conjecture is known as Price's Law. There has been much mathematical interest in studying pointwise decay rates of waves on relativistic geometric backgrounds, including the Schwarzschild and Kerr spacetimes (the latter describes the geometry resulting from a rotating black hole). Both these geometries tend toward flat at a rate of $r^{-1}$. 

Price's Law was proved in \cite{Ta:13} (see also \cite{DoScSo:11} for
the mode-by-mode estimate as well as \cite{DaRo:10} and \cite{Lu:12}
for earlier decay estimates). Pointwise decay rates for the Kerr
spacetime were studied in \cite{dafrod} and \cite{finster}. The
techniques in \cite{DoScSo:11}, \cite{Ta:13}, \cite{Mo:20}, and the
current work involve taking the Fourier transform in time and
therefore do not readily extend to non-stationary geometries. In
\cite{metato} the authors proved Price's Law for non-stationary
asymptotically flat spacetimes and established the $t^{-3}$ decay rate
for a class of perturbations of the Kerr spacetime.  Results similar
to the current work but also considering non-stationary spacetimes
were obtained concurrently by Looi in \cite{looi}. The author
  uses vector field techniques and an iterative argument to obtain
  decay rates which depend on the rate at which the background
  geometries become flat.  While the results of \cite{looi} are in
  fact stronger than those obtained here in the sense that they apply
  to nonstationary metrics (and do not lose an $\epsilon$ power), the
  results of this paper give direct and precise estimates in the
  frequency domain---in particular, near zero frequency---which we
  hope will find further applications.  For instance, given a
  situation in which high-frequency resolvent estimates are known
  directly (and there are many such situations, often with estimates
  obtained directly in frequency domain by semiclassical propagation
  of singularities; see e.g.\ \cite{Wu:12}), we can easily combine the
  high- and low-frequency estimates to obtain new decay estimates in
  the time domain.

In the present work as well as in \cite{Mo:20} and \cite{Ta:13} an
integrated local energy decay estimate is assumed to hold. Such
estimates were established on the Schwarzschild geometry in
\cite{blue2003}, \cite{dafermos2009red}, and \cite{MaMeTaTo}. For the
Kerr spacetime with low angular momentum, local energy estimates were
proved in \cite{andblue15}, \cite{dafermos2009red}, and
\cite{dafrodshlap16}; the major challenge here was the trapping of
null geodesics within a compact set in the space variables.  The
assumptions in \cite{Ta:13} therefore hold for Schwarzschild and Kerr
with low angular momentum.  We discuss the integrated local energy
decay estimates in more detail later in the introduction. 

Finer analysis of the asymptotics of solutions on Schwarzschild space
was obtained in \cite{AnArGa:18}, including a characterization of when the
$t^{-3}$ decay rate is a \emph{lower} bound for the decay.

A different approach to Price's law was pioneered by Hintz
\cite{Hi:20}, who showed that the estimate in Price's Law on Kerr
backgrounds is 
sharp, and obtained explicit 
leading order asymptotics.  Here, rather than use vector field
  methods in physical space, the author employs resolvent
estimates, after Fourier transforming the equation in the appropriate variables.
No explicit local energy decay assumption is employed: the
author instead considers metrics for which the associated spectral family (given by formal Fourier transform in time) satisfies appropriate hypotheses on its inverse (the resolvent). Integrated local energy decay estimates
are in fact intimately related to resolvent estimates. For example, in \cite{MeStTa} the authors establish a full
spectral characterization of local energy decay in the context of
nontrapping asymptotically flat spacetimes.

An analysis of the asymptotics of solutions that distinguishes the
contributions due to low angular modes from the more decaying
contributions from higher ones
has recently been carried out on
Reissner--Nordstr\"om backgrounds in \cite{AnArGa:21a} and on Kerr
backgrounds in \cite{AnArGa:21b}.

Our approach blends the hypotheses of \cite{Ta:13} and associated
works with an adaptation of the new approach of Hintz \cite{Hi:20}, which in turn
harnesses powerful low frequency resolvent estimates recently developed by
Vasy \cite{Va:20}. The high frequency estimates thus follow
\cite{Ta:13} and \cite{Mo:20}, with the energy decay assumption being
the crucial hypothesis ensuring that in return for enough derivatives
of regularity, this part of the solution in fact decays at any desired
rate.  At low frequency, by contrast, the local energy decay
assumption yields absence of resonances in the upper half-space, which
is certainly one necessary condition for decay estimates to hold (and
this is an explicit spectral hypothesis, e.g.\ in \cite{Hi:20}).
We obtain additional information about the low-frequency asymptotics of
asymptotically Euclidean resolvents, which is essentially
independent of the small-scale geometry, via estimates from
\cite{Hi:20} and \cite{Va:20}: we apply these iteratively to obtain conormal
estimates for the resolvent at zero frequency.

A crucial step in our analysis is thus writing an expansion of the
resolvent at zero frequency (see Lemma \ref{lemma:P0inverse}). Low frequency resolvent expansions of the
Laplacian were first studied in \cite{JeKa}. The geometric context we
consider here reduces to analyzing perturbations of the flat Laplacian
where the perturbation depends on the spectral parameter. The presence
of such terms arises from the fully
Lorentzian nature of the perturbations considered here. Closely related results on low frequency
spectral behavior and local wave decay were previously obtained in the setting of scattering manifolds in
\cite{GuHaSi:13}, where the authors use an extremely precise description of the spectral measure
at low frequencies to establish decay rates dependent on the dimension
of the spacetime.

Questions similar to those treated here were studied in
\cite{bonyhaf2} and \cite{bonyhaf1} where the authors established
local decay rates for waves on asymptotically flat, stationary
spacetimes which tend toward flat at different rates. There are
several key differences compared with the current work. First, as
noted above, we handle full Lorentzian perturbations of flat Minkowski
space rather than restricting to perturbations of the Laplacian. This
leads to the metrics considered in this paper containing $dtdx_i$
terms, which results in mixed space-time differential operators in our
wave operator. Second, we allow for the possibility of unstable
trapping on our background. In \cite{bonyhaf2} and \cite{bonyhaf1}, a
nontrapping assumption is used in order to obtain decay for the high
frequency part of a solution to the wave equation (it is not needed
for the low frequency part). Third, our result improves upon the
established decay rates. Finally we note that \cite{bonyhaf1}
considers $(1+n)$ dimensional geometries for $n \ge 2$ and
\cite{bonyhaf2} considers $n$ odd with $n \ge 3$. The current work
only studies $(1+3)$ dimensional spacetimes.

\begin{subsection}{Wave Equation}
	On flat Minkowski space, the wave operator, denoted $\Box$, is
        given by $\Box = \partial_t^2 - \Delta = \partial_t^2 -
        \sum_{i=1}^3 \partial_{x_i}^2$. More generally, on a
        Lorentzian spacetime with metric $g = g_{\alpha\beta}$, the
        wave operator is denoted $\Box_\gmetric$ and is given by the
        d'Alembertian 
	$$
		\Box_\gmetric = \frac{1}{\sqrt|\gmetric|} \partial_\alpha \sqrt{|\gmetric|}\gmetric^{\alpha \beta} \partial_\beta
	$$
	where $g^{\alpha\beta}$ are the dual metric coordinates. 
	\end{subsection}

\begin{subsection}{Local Energy Decay} \label{subsection:LED}
	On stationary spacetimes, as considered here, solutions to the wave equation may have constant energy. If the background geometry allows solutions to spread out, then energy may decay within compact sets. We assume an integrated local energy estimate where the local decay is fast enough to be integrable in time. Estimates on the local decay of energy have a long history dating back to the work of Morawetz \cite{mora}. The specific version of the estimates used here come from the work of Metcalfe and Tataru in \cite{MeTa:12}. In addition to local energy decay, we will assume uniform energy bounds on solutions of the wave equation. Such bounds necessarily hold on stationary spacetimes where the Killing vector field $\partial_t$ is everywhere time-like. The uniform energy estimate also holds on Schwarzschild and Kerr, although in these geometries $\partial_t$ is not everywhere time-like (see e.g. \cite{dafermos2009red}, \cite{MaMeTaTo}, \cite{DafRod11} among others).
	
	We consider the Cauchy problem
	\begin{equation} \label{inhcauchy}
		(\boxg + V)u = f, \qquad u(0,x) = u_0, \qquad \partial_t u(0,x) = u_1	
	\end{equation}
where $V$ is a suitably decaying potential (see the statement of the main theorem for conditions on $V$). The Cauchy data at time $t$ is denoted $u[t] = \Big(u(t,\cdot), \partial_tu(t,\cdot)\Big)$.
	
		\begin{definition}	
			We say the evolution \eqref{inhcauchy} satisfies the uniform energy bounds if:
				\begin{equation} \label{unifenbd}
					\| u[t]\|_{\dot H^{k,1} \times H^k} \leq c_k (\|u[0]\|_{\dot H^{k,1} \times H^k} +\|f\|_{ L^1H^k}),
	\qquad t \geq 0, \quad k \geq 0.
				\end{equation}
		\end{definition}	
Here $H^k$ denotes the usual Sobolev space, and we say $\phi \in \dot{H}^{k,1}$ if $\nabla \phi \in H^k$, and write
$$
\norm{\phi}_{\dot{H}^{k,1}}= \norm{\nabla \phi}_{H^k}.
$$
	
In the following definitions we use $\partial$ to denote the space-time gradient while $\nabla$ is reserved for the gradient in spatial variables only. We write $\langle r \rangle := \sqrt{1+r^2}$ and define the dyadic region $A_m := \{ x : 2^m \le \langle r \rangle \le 2^{m+1}\}$.  The local energy norm we use is defined by
			\[
				\|u\|_{LE} = \LEnorm{u}.
			\]
		Its $H^1$ analogue is given by
			\[
				\|u\|_{LE^1} = \LEonorm{u},
			\]
		and the dual norm is given by
			\[
				\|f\|_{LE^*} = \LEsnorm{f}.
			\]
		For functions with higher regularity we define the following norms
			\[
				\|u\|_{LE^{1,N}} = \sum_{j \le N} \|\partial^j u\|_{LE^1}, \qquad \|f\|_{LE^{*,N}} = \sum_{j \le N} \|\partial^j f\|_{LE^*}.
			\]
		The spatial counterparts of the $LE$ and $LE^*$
                space-time norms are
			\[
				\|v\|_\LE = \sup_m \| \la r \ra^{-\frac{1}{2}} v\|_{L^2(A_m)}; \qquad \|g\|_\LEs = \sum_m \| \la r \ra^{\frac{1}{2}} g\|_{L^2(A_m)}
			\]
		with the higher regularity norms defined by
			\[
				\|v\|_{\mathcal{LE}^{N}} = \sum_{j \le N} \|\nabla^jv\|_{\mathcal{LE}}, \qquad \|g\|_{\mathcal{LE}^{*,N}} = \sum_{j \le N} \|\nabla^jg\|_\LEs.
			\]			

		\begin{definition}
			We say the evolution \eqref{inhcauchy} satisfies the local energy decay estimate if:
				\begin{equation} \label{led}
	 				\| u\|_{LE^{1,N}} \leq c_N (\|u[0]\|_{H^{N,1} \times H^N} + \|f\|_{LE^{*,N}}), \qquad N \geq 0.
				\end{equation}
		\end{definition}
	
                Local energy decay is known to hold in several
                nontrapping geometries. For sufficiently small
                perturbations of flat space without trapping, local
                energy decay was established in \cite{alinhac},
                \cite{metsog06}, and \cite{mettat1}. The case of
                stationary product manifolds was considered in
                \cite{burq1998}, \cite{bonyhaf1}, and
                \cite{sogwan}. The nontrapping case was studied more
                generally in \cite{MeStTa}. If trapping occurs then
                the local energy decay estimate does not hold
                (\cite{ralston},
                \cite{sbierski2015characterisation}). However, if the
                trapping is sufficiently unstable (i.e.\, perturbing
                a trapped geodesic typically results in geodesics which
                escape to infinity) then a weaker form of local energy
                decay may hold (see, e.g., \cite{WuZw:11} in the case
                of the normally hyperbolic trapping that occurs in
                Kerr black hole spacetimes). In the case of trapping, there is
                necessarily a loss of derivatives on the right hand
                side of the estimate (see e.g. \cite{bcmp}).
			\begin{definition}
	 		We say the evolution \eqref{inhcauchy} satisfies the weak local energy decay estimate if:
				\begin{equation} \label{wled}
	 				\| u\|_{LE^{1,N}} \leq c_N (\|u[0]\|_{\dot H^{N+\loss,1} \times H^{N+\loss}} + \|f\|_{LE^{*,N+\loss}} ), \qquad N \geq 0
				\end{equation}
	for some $\loss > 0$.
		\end{definition}

This weak local energy decay estimate is generally obtained by using a cutoff function to remove the trapped set. The precise derivative loss depends on the trapping. Here we allow for an arbitrary but fixed loss, $\loss$.
\end{subsection}

\begin{subsection}{b-Sobolev spaces}
Our hypotheses on the Cauchy data and our key estimates in the iteration at
low frequency will be stated in the language of \emph{b-Sobolev
  spaces}.  These spaces, described in detail in \cite{Melrose:APS}
(and discussed further in Section~\ref{section:background} below),
can be easily described in the context at hand, at least for positive integer
orders, as the spaces of functions enjoying Kohn--Nirenberg symbol
estimates to finite order with respect to
$L^2$ (rather than the usual $L^\infty$):
$$
u \in H_b^{m,l}(\overline{\RR^3}) \Longleftrightarrow
\smallnorm{\ang{x}^{l+\smallabs{\alpha}} \pa_x^\alpha u}_{L^2(\RR^3)} \leq C_\alpha,\quad
\abs{\alpha}\leq m.
$$
The fact that the space is nominally defined on $\overline{\RR^3},$
the radial compactification of $\RR^3,$ is a nod to the fact that
these spaces are more generally defined on manifolds with boundary;
details follow in Section~\ref{section:background} and
Appendix~\ref{app:bgeometry} below.

We will also use in our hypotheses the usual spaces of
Kohn--Nirenberg symbols, defined by the estimates
\begin{equation}\label{KN}
a\in S^l(\RR^n) \Longleftrightarrow a\in \CI(\RR^n),\ \sup \ang{x}^{-l+\abs{\alpha}}
  \abs{\pa_x^\alpha a} <C_\alpha\quad \text{ for all } \alpha.
\end{equation}
\end{subsection}

\begin{subsection}{Statement of Main Theorem}

  Let $\gmetric,$ $\hmetric$ satisfy the hypotheses of Section 1.5 of \cite{Mo:20}:

	\begin{enumerate}
		\item[(i)] $\gmetric$ is stationary (i.e.\ the metric coefficients are time independent).
		\item[(ii)] The submanifolds $t = constant$ are space-like (i.e.\ the induced metric on the spatial submanifolds is positive definite).
		\item[(iii)] Let $\kappa \in (1,\infty)\setminus \N$. The metric $\gmetric$ is asymptotically flat in the sense that $\mathfrak{g}$ can be written as
	\[ \gmetric = \mathfrak{m} + \hmetric \]
where
	\[\hmetric = \hmetric_{00}(x)dt^2 + \hmetric_{0i}(x)dtdx_i +  \hmetric_{ij}(x)dx_idx_j \]
	with $\hmetric_{\alpha\beta} \in S^{-\kappa}(\RR^3)$ for $\alpha, \beta \in \{0,1,2,3\}$.
	\end{enumerate}

\begin{theorem}
	Let $u$ solve the homogeneous Cauchy problem
	$$
		(\Box_g + V)u = 0, \quad u(0,x) = u_0, \quad \partial_t u(0,x) = u_1
	$$
	for $V \in S^{-\kappa-2}(\R^3)$ and
        $$
u_0 \in H_b^{s+1, \kappa+7/2},\quad u_1 \in H_b^{s, \kappa+7/2},
$$
with
  $$
  s> (\loss+1)(2\kappa + 9) + 2.
  $$
         Assume the evolution \eqref{inhcauchy} satisfies the uniform
         energy bounds \eqref{unifenbd} and the weak local energy
         decay estimate \eqref{wled}. Then $|u(t,x)| \le C_\ep
         t^{-\kappa-2+\epsilon}$ for any $\epsilon >0$, uniformly on
         compact sets in $x$.
\end{theorem}
\end{subsection}

\begin{subsection}{Paper Overview}
The operator $\Box_\gmetric +V$ can be replaced by an operator of the form
	$$ 
		P = \partial_t^2 - \Delta + \partial_tP^1+P^2
	$$
where 
	$$
		P^1 \in S^{-\kappa}\partial_x + S^{-\kappa-1}  \quad \hbox{and} \quad P^2 \in  S^{-\kappa} \partial_x^2 + S^{-\kappa-1}\partial_x + S^{-\kappa-2}.
	$$
This is obtained by working in normalized coordinates and multiplying by $(\gmetric^{00})^{-1}$ so the coefficient on $\partial_t^2$ is 1. We refer the reader to \cite{Mo:20} for details of the calculation. We work with this operator throughout. The resolvent, $P_\sigma^{-1}$, associated to $P$ is given formally by inverse Fourier transforming the operator $P$ and taking the inverse. 

In section \ref{EnergyToSpectralInfo} we define the resolvent and
extract relevant spectral information from the energy assumptions. In
particular, we show the resolvent is well defined in the upper half
plane and extends continuously to the real axis. The results in section \ref{EnergyToSpectralInfo} allow us to relate the inverse Fourier transform of a solution $u$ to the initial data of the homogeneous Cauchy problem via the resolvent:
	$$
		\check{u}(\sigma, x) = (2\pi)^{-\frac{1}{2}}P_\sigma^{-1}(-i\sigma u_0 +P^1u_0-u_1)
	$$
and to recover $u$ by taking the Fourier transform (with integration along $\sigma \in \R$). The amount of decay we are able to obtain then depends on conormal regularity estimates (i.e.\ bounds on $(\sigma \partial_\sigma)^j P_\sigma^{-1}g$ for appropriately chosen $g$). We handle the high and low frequency cases separately. The high frequency part of $u$ is sensitive to the trapping dynamics, which are controlled by our weak local energy decay assumption. Indeed, the spectral information derived from the energy assumptions is sufficient to handle the high frequency part of the solution to the wave equation, and we find that in exchange for enough derivatives on the Cauchy data, we could obtain any desired polynomial time-decay for this piece of the solution. The low frequency part of $u$, by contrast, is sensitive to the far away behavior of the background geometry, and it is this latter piece which ultimately dictates the final decay rate.

We establish conormal
estimates for the high energy resolvent in section
\ref{HighEnergy}. We then turn to the low energy analysis in section
\ref{LowEnergy}. This includes deriving an expansion of the resolvent
at zero energy, which is then used to find a helpful expression for the
resolvent at low frequencies. The low frequency analysis utilizes conormal and b-Sobolev spaces, and we provide an overview of these function spaces in section \ref{section:background}. Finally, we prove the main theorem in section
\ref{MainProof}.

\end{subsection}

\subsection{Acknowledgements} The authors are grateful to Peter Hintz
for numerous helpful conversations, as well as to Jason Metcalfe,
Mihai Tohaneanu, and Shi-Zhuo Looi.  Two anonymous referees also
supplied helpful comments and corrections.

The second author gratefully acknowledges support from Simons Foundation grant 631302. The first author gratefully acknowledges the support of NSF Postdoctoral Fellowship DMS-2002132.

\section{Background on Function Spaces}\label{section:background}

In dealing with low-frequency estimates, it will be convenient to
compactify our asymptotically Euclidean space and to
employ the language of conormal and b-Sobolev spaces on manifolds with
boundary, as introduced by Melrose--Mendoza
\cite{MeMe:83} based on Melrose's b-calculus of pseudodifferential
operators \cite{Me:81}. (See \cite{Melrose:APS} for an extended exposition.)
Some details of the local characterizations of these spaces near the
boundary at infinity via Mellin transform have been relegated to Appendix~\ref{app:bgeometry}.

Let $\overline{\RR^n}\simeq \overline{B^n}$ denote the radial compactification of Euclidean
space to the unit ball, with $\RC:\RR^n \to
(\overline{\RR^n})^\circ$ given by the compactification map
$$
\RC(x) =\frac{x}{1+\smallang{x}}.
$$
Note that $\RC^*(\CI(\Rnbar))=S^0_{\text{cl}}(\RR^n),$ the space of
``classical'' symbols on $\RR^n$ (those satisfying \eqref{KN} with
$l=0$ and additionally admitting an asymptotic
expansion in negative integer powers of $\smallabs{x}$).  The function
$$\rho \equiv \smallabs{x\circ\RC^{-1}}^{-1}$$ extends (except for a
singularity at the origin) to a 
smooth function on $\overline{\RR^3}$ that is a boundary defining
function (it vanishes to exactly first order at the boundary).  We
will freely employ the abuse of notation
$$
\rho=r^{-1}=\smallabs{x}^{-1},
$$
ignoring the pushforward/pullback by $\RC,$ and moreover will consider $\rho$ to
be extended to a globally smooth function on $\overline{\RR^3},$
eliminating the singularity at $x=0.$
   
More generally, we temporarily let $X$ denote any manifold with
boundary.  Let $\V_b(X)$ denote the
   space of ``b-vector fields,'' i.e.\, those which are tangent to $\pa
   X.$  If $(\rho,y)$ are coordinates in a collar neighborhood of $\pa
   X,$ with $\rho$ a boundary defining function and $y$ coordinates on
   $\pa X,$ extended to the interior, then $\V_b(X)$ is locally
   spanned by $\rho \pa_\rho, \pa_y$ over $\CI(X);$ in the special
   case $X=\Rnbar$ these vector fields correspond to $-r \pa_r,
   \pa_\theta;$ note in particular that their norm is $O(r)$ as
   $r \to\infty.$

The b-differential operators, $\hbox{Diff}_b^m(X)$, are defined as the
$\mathcal{C}^\infty$-span of up to $m$-fold products of vector fields in $\mathcal{V}_b(X).$
Given a fixed volume form on $X^\circ$ (possibly singular at $\pa X$) we define $L^2(X)$  with respect to
the volume form, and then set, for $k \in \NN,$
$$
H_b^k(X) = \big\{u \in L^2(X)\colon V_1\dots V_j (u) \in L^2(X)\ \forall
j\leq k, V_i \in \V_b(X)\big\}.
$$
We can further define $H_b^s(X)$ for $s \in \RR$ by interpolation and
duality, or by use of the calculus of b-pseudodifferential operators
microlocalizing the algebra $\Diffb(X)$ as in \cite{Melrose:APS}.

\emph{In this paper, we will always employ the standard metric volume
  form on $\RR^3$} in defining $L^2$ and Sobolev spaces on
$\overline{\RR^3}$.  This volume form is given, near the boundary $\rho=0$, by
$$
dV=\frac{d\rho \, d\angvar}{\rho^4}
$$
(with $d\angvar$ shorthand for the volume
for on $S^2$) hence appears quite singular on the
compactification.\footnote{We caution the reader that this choice
of convention is not universal in the subject, with the ``b-volume
form'' $d\rho\,  d\angvar/\rho$ also having a considerable popularity,
since its use would eliminate the factors of $\rho^{3/2}$ that
bedevil the accounting used here.}  Note that Sobolev
spaces with no subscripts will continue to denote ``ordinary'' Sobolev
spaces on $\RR^3.$  These can, if desired, be identified with the
``scattering Sobolev spaces'' on $\overline{\RR^3}$ as introduced by
Melrose \cite{Me:94}, but we will not employ this terminology below.

We generalize the b-Sobolev spaces to \emph{weighted} b-Sobolev spaces
by simply setting
$$
H_b^{m,\ell}(X) = \rho^\ell H_b^m(X).
$$

A related function space to the b-Sobolev space is that of
\emph{conormal functions} enjoying infinite-order iterated regularity
under $\mathcal{V}_b(X),$ measured with respect to $L^\infty$ rather
than $L^2:$
$$
u \in \A^\gamma \Longleftrightarrow V_1\dots V_N (u) \in \rho^\gamma L^\infty(X)
\text{ for all } N \in \NN,\ V_j \in \mathcal{V}_b(X).
$$
Such estimates are closely related to the standard Kohn--Nirenberg
symbol estimates, since $\mathcal{V}_b(\overline{\RR^3})$ is spanned
over $\CI$ by $r\pa_r=-\rho\pa_\rho,$ and $\pa_\theta$; thus
radial compactification gives an
isomorphism
$$
\A^\gamma (\overline{\RR^3})\equiv S^{-\gamma}(\RR^3),
$$
with the symbol spaces $S^\bullet$ defined by \eqref{KN} above.
The point is that the definition of the symbol spaces can be rephrased as iterated regularity under
the vector fields
$\smallang{x} \pa_{x_j},$ which span $\mathcal{V}_b.$

The conormal spaces are also closely related to $H^{\infty,l}_b(X),$ but with a vexing
shift in orders owing to the metric volume form:
\begin{equation}\label{sobemb}
		H_b^{\infty,\ell}(X) \subset \A^{\ell+3/2}(X) \subset
                H_b^{\infty,\ell-}(X),\quad \ell \in \RR.
              \end{equation}

\section{Spectral Information from Energy Assumptions} \label{EnergyToSpectralInfo}

We define a spectral family associated to $P$ by
\begin{equation}\label{Psigma}
		\begin{split}
		P_\sigma &\equiv e^{it\sigma}Pe^{-it\sigma}\\
			 &= -\sigma^2 - \Delta - i\sigma P^1 + P^2
		\end{split}
\end{equation}
                which is equivalent to formally inverse Fourier transforming the
operator in $t$. In this section we use the energy assumptions to study the existence and boundedness of the operator $P_\sigma^{-1}$.

We will use the unitary normalization for the Fourier transform on $\RR^n:$
$$
(\F f)(\xi)=(2\pi)^{-n/2} \int f(x) e^{-i \xi x}\, dx,
$$
with the inverse Fourier transform then given by its formal adjoint
(and denoted $\F^{-1} f$ or $\check{f}$).

\begin{proposition} \label{prop:resex} If $\Im \sigma > 0$, the
  operator $P_{\sigma}: H^2 \to L^2$ is invertible.  Furthermore, if $u$ satisfies
  \eqref{inhcauchy} then for $\Im \sigma >0$ we have
		\begin{equation} \label{hatures}
			(2\pi)^{1/2}\check{u}(\sigma,x) = 
                        P_{\sigma}^{-1} ((2\pi)^{1/2}\check{f}(\sigma) -i\sigma u_0 + P^1 u_0 - u_1).
		\end{equation}
	\end{proposition}
	
	\begin{proof}
Let $u$ solve $Pu=f$ for $t\geq 0.$  By the uniform energy bound \eqref{unifenbd},
$e^{it \sigma} u$ enjoys exponential energy decay at $t \to +\infty$
for $\Im \sigma>0,$ hence, with $H(t)$ denoting the Heaviside function,
$\F_{t\to \sigma}^{-1} (H(t) u)$ is an analytic function of $\sigma$
in the upper half-space, taking values in the energy space.

          Integrating by parts in $t$ moreover gives, for $\Im \sigma>0,$
\begin{align*}
(2\pi)^{1/2}\F^{-1} ( H(t) &Pu)(\sigma,x)= \int_0^\infty e^{+it
                                          \sigma} Pu\, dt\\ &= (2\pi)^{1/2}P_\sigma \F^{-1} (H(t)u)(\sigma) -\pa_t u(0) +i\sigma u(0) - P_1 u(0)
\end{align*}
Consequently, 
\begin{align*}
(2\pi)^{1/2}P_\sigma \F^{-1} (&H(t)u)(\sigma)\\
& =-(2\pi)^{1/2}\F^{-1} ( H(t) Pu)(\sigma,x) +\pa_t u(0) -i\sigma u(0) + P_1 u(0).
\end{align*}
Thus we can solve $P_\sigma v=g,$ for $g \in L^2,$ by solving the IVP
\eqref{inhcauchy} with
$u_0=0,$ $u_1=g$
and setting
\begin{equation}\label{vdef}
v(t,x) = \int_0^\infty u(t,x) e^{it\sigma} \, dt.
  \end{equation}
(A priori this construction
would produce $v \in H^1,$ but $v \in H^2$ then follows by ellipticity
of $P_\sigma$.)

Now we turn to injectivity of $P_\sigma.$  If $P_\sigma w=0$ with $\Im
\sigma>0,$ we set
$$
u(t,x) = e^{-it \sigma} w(x).
$$
Then by \eqref{Psigma} we obtain
$$
P u = e^{-it\sigma} P_\sigma w=0.
$$
Since $\Im
\sigma>0,$ $\norm{u}$ is exponentially growing as $t\to +\infty,$
contradicting the uniform energy bounds.  Hence we obtain injectivity
of $P_\sigma.$	\end{proof}

The $\mathcal{LE}_{\sigma}$ norm, in which we measure the resolvent $v = P_\sigma^{-1} g$, is defined by
	\begin{equation}
		\|v\|_{\mathcal{LE}_{\sigma}^N} = \| (|\sigma| + \la r \ra^{-1})v \|_{\mathcal{LE}^{N}} + \| \nabla v \|_{\mathcal{LE}^{N}} + \|(|\sigma| + \la r \ra^{-1})^{-1}\nabla^2 v \|_{\mathcal{LE}^{N}}. \label{letnorm}
	\end{equation}
Next we transfer bounds in the local energy decay estimate to bounds on the resolvent measured in the $\mathcal{LE}_\sigma$ norm which hold down to the real axis. Ultimately we wish to obtain decay rates for $u$ by taking the Fourier transform of \eqref{hatures}, but integrating in $\Im \sigma >0$ would lead to exponential blow-up in time. The extension of $P_\sigma^{-1}$ to $\sigma \in \R$ allows us to integrate along $\sigma \in \R$. 

The following proposition is analogous to results in \cite{Ta:13} (see Proposition 9 and Corollary 12) where $\loss = 3$ (see equation (4.2)). We provide an outline of the proof, focusing on how our assumed $\varpi$ derivative loss in the weak local energy decay estimate \eqref{wled} affects the derivative loss in the estimate \eqref{resest}.
	
	\begin{proposition} \label{Prop:resest}
		If $\Im \sigma \ge 0$ and $g \in \mathcal{LE}^{*,N+\loss+1}$ for fixed $N \in \N$, then $v = P_\sigma^{-1} g$ satisfies
		\begin{equation} \label{resest}
			\|v\|_{\mathcal{LE}_\sigma^N} \lesssim \|g\|_{\mathcal{LE}^{*,N+\loss+1}}.
		\end{equation}
	\end{proposition}
	
	\begin{proof}
Recall that we may construct $v$ for $\Im \sigma>0$ by solving the IVP
with Cauchy data $(0,g)$ and evaluating the integral \eqref{vdef};
this construction continues to make sense distributionally down to
$\Im \sigma=0,$ indeed with explicit weighted Sobolev estimates, as we will now see.

          The weak local energy decay estimate allows us to establish the inequality
		\begin{equation} \label{wledtrans2}
			\| (\la r \ra^{-1} + |\sigma|) v \|\len{N} + \| \nabla v \|\len{N} \lesssim \sum_{j \le N+\loss} (1+|\sigma|)^{N+\loss-j}\| \nabla^j g \|\les, 
		\end{equation}
Note \eqref{wledtrans2} follows formally from \eqref{wled} by Plancherel after inverting the Fourier transform. From \eqref{wledtrans2} we are able to bound the first two terms in the $\mathcal{LE}_\sigma^N$ norm:
	 \begin{equation} \label{resbnda}
	 	\| ( \la r \ra^{-1}+|\sigma| )v \|_{\mathcal{LE}^{N}} + \| \nabla v \|_{\mathcal{LE}^{N}} \lesssim \|g\|_{\mathcal{LE}^{*,N+\loss}}.
	 \end{equation} 
We refer the reader to \cite{Ta:13} for details of the process (cf.\ equation (4.6)). 
		
		Now consider the second order term in the $\mathcal{LE}_\sigma^N$ norm. If $(\la r \ra^{-1} +|\sigma|)^{-1}\lesssim 1$ then
			\[ \| (\la r \ra^{-1} + |\sigma|)^{-1} \nabla^2 v \|\len{N} \lesssim \| \nabla^2 v \|\len{N} \lesssim \|\nabla v\|\len{N+1} \lesssim \| g \|\lesn{N+\loss+1}. \]
This step is where the extra derivative loss in \eqref{resest} versus \eqref{wled} comes from.
		It is left to consider the case where $\la r \ra$ is large and $|\sigma| \lesssim 1$. Here we write
			\[ \| (\la r \ra^{-1} +|\sigma|)^{-1} (-\Delta + P^2)v \|\len{N} \lesssim \| |\sigma| v \|\len{N} + \| P^1 v \|\len{N} + \|\la r \ra g\|\len{N}. \]
		The first two terms have already been shown to satisfy the desired bounds. For the third term, straightforward calculation yields $\| \la r \ra g \|\len{N} \lesssim \| g \|\lesn{N}$. The estimate is transferred to $\nabla^2 v$ using standard elliptic arguments.
	\end{proof}

\section{High Energy Conormal Estimates} \label{HighEnergy}
The goal of this section is to obtain pointwise bounds on $(\sigma\partial_\sigma)^MP_\sigma^{-1}g$ for $g \in \mathcal{LE}^*$. The results in this section are direct analogues to results in \cite{Ta:13} and \cite{Mo:20}. We provide sketches of the arguments, but detail how the regularity requirements depend on the loss, $\varpi$, in the weak local energy decay estimate \eqref{wled}.  

We consider the vector fields $T \in \{ \partial_{x_i} | i = 1,2,3\},
\Omega \in \{ x_i\partial_{x_j} - x_j\partial_{x_i} | i,j = 1,2,3 \}$,
and $S = r\partial_r - \sigma \partial_\sigma$. Note if $g$ is
independent of $\sigma$ then $Sg = r\partial_r g$.  Since
$r\pa_r=-\rho\pa_\rho,$ this implies that
$$
T^i\Omega^jS^k g \in L^2 \text{ for all } i+j+k \leq M
\Longleftrightarrow g \in \Hb^{M,0}.
$$
(Note that we include the $T$ derivatives, which are powers of vector
fields in $\rho \mathcal{V}_b,$ merely to ensure differentiability at $r=0.$)
The assumptions here are stated in terms of the energy
space $\mathcal{LE}^*$, and we note $\mathcal{LE}^* \subset
r^{-\frac{1}{2}}L^2$.

In Lemmas \ref{Lemma:ijkresest} and \ref{Lemma:vijkptbnd} we take $Q$ to represent an operator of the same form as $\sigma P^1 + P^2$ but let the exact coefficients change each time $Q$ appears. 

The following lemma is an analogue to \cite{Ta:13} Proposition 10.

	\begin{lemma}\label{Lemma:ijkresest}
		If $\Im \sigma \ge 0$ and $g \in \mathcal{LE}^*$ satisfies
		\begin{equation} \label{gijkbnd}
			\|T^i\Omega^jS^k g\|_{\mathcal{LE}^*} \lesssim 1, \quad i+ (\loss+1)j + (\loss +1)k \le M
		\end{equation}
		for some positive integer $M$, then $v=P_\sigma^{-1}v$ satisfies
		\begin{equation}\label{vijkbnd99}
			\|T^i\Omega^jS^k v\|_{\mathcal{LE}_\sigma} \lesssim 1, \quad i+ (\loss+1)j + (\loss +1)k \le M - \loss -1
		\end{equation}
	\end{lemma}
	
	\begin{proof}
		Applying Proposition \ref{Prop:resest} to
			$$ P_\sigma \vijk = \gijk + [P_\sigma, \ijk]v $$
		yields
			$$ \| \vijk \|\letn{N} \lesssim \| \gijk \|\lesn{N + \loss + 1} + \| [P_\sigma,\ijk] v \|\lesn{N + \loss + 1}. $$
		Since $\|\vijk\|\letn{} = \|\Omega^jS^kv\|\letn{i}$ we only require bounds on $\|\Omega^j S^k v\|\letn{N}$. We illustrate the general method with concrete examples and highlight the role of $\loss$ in determining the requisite regularity.
		
		Consider $\Omega v$. Direct calculation shows
                $[P_\sigma, \Omega] = Q$ and
		$$
			\| Q \phi \|\lesn{N} \lesssim \| (\langle r \rangle^{-1} + |\sigma|) \phi \|\len{N} + \|\nabla \phi \|\len{N} =:  \| \phi \|\lettn{N}
		$$
		for $\phi$ with sufficient regularity and decay. Note that
                the $\mathcal{LE}_{\tau,1}^N$ norm defined here is the
                first two terms of the $\mathcal{LE}_\tau^N$ norm \eqref{letnorm},
                and
                by \eqref{resbnda} we have $\|v\|\lettn{N} \lesssim
                \|g\|\lesn{N+ \loss}$. Thus we find 
		\begin{align*}
			\|\Omega v \|\letn{N} &\lesssim \| \Omega g \|\lesn{N+\loss+1} + \|Q v\|\lesn{N+\loss+1} \\
				& \lesssim 1 + \| v \|\lettn{N + \loss + 2} \\
				& \lesssim 1 + \|g \|\lesn{N+2(\loss+1)}.
		\end{align*}
	We see there is one loss of $\loss +1$ due to the estimate in Proposition \ref{Prop:resest} and a subsequent loss of $\loss +1$ for each $\Omega$. This justifies the requirement $i + (\loss +1)j \le M-\loss - 1$.
	
	Next we consider $S v$. Note $[P_\sigma, S] = 2P_\sigma + Q$. We calculate
	\begin{align*}
		\| SP_\sigma^{-1}g\|\letn{N} &\lesssim \|g\|\lesn{N+\loss+1} + \|Sg\|\lesn{N+\loss+1} + \|Q v\|\lesn{N+\loss+1} \\
			&\lesssim 1 + \|v\|\lettn{N+\loss+2}\\
			&\lesssim 1 + \|g\|\lesn{N + 2(\loss+1)}.
	\end{align*}
As above, this demonstrates the requirement $i + (\loss+1)j + (\loss +1)k \le M- \loss -1$. The general case follows by induction. Details are omitted here. 
	\end{proof}

 In Lemma \ref{Lemma:vijkptbnd} and Proposition \ref{largesigmaresbnd} we will make use of the following result, which we quote from \cite{Ta:13} Proposition 11:

\begin{lemma}\label{Lemma:radcond}
	If $\sigma \in \mathbb{R} \setminus \{ 0\}$ and $g \in \mathcal{LE}^{*,4}$, then $v = P_\sigma^{-1}g$ satisfies the outgoing radiation condition
	\begin{equation}\label{radcondsimp}
		\lim_{m \to \infty} 2^{-\frac{m}{2}}\| (\pa_r-i\sigma)v\|_{L^2(A_m)} = 0.
	\end{equation}
\end{lemma}

A consequence of Lemma \ref{Lemma:radcond} and the proof of Lemma \ref{Lemma:ijkresest} is that for $g$ satisfying \eqref{gijkbnd}, the radiation condition holds for $T^i\Omega^jS^kv$ with appropriate values of $i, j,$ and $k$:
\begin{equation}\label{radcond}
	\lim_{m \to \infty} 2^{-\frac{m}{2}}\| (\pa_r-i\sigma)T^i\Omega^jS^kv\|_{L^2(A_m)} = 0, \quad i+(\loss+1)j + (\loss+1)k \le M - 4.
\end{equation}
As in the proof of Lemma \ref{Lemma:ijkresest}, we see there is one loss of $\loss +1$ due to Lemma \ref{Lemma:radcond} and a subsequent loss of $\loss +1$ for each $\Omega$ and each $S$ due to the commutator terms $[P_\sigma, \Omega]$ and $[P_\sigma, S]$.

Now we provide preliminary pointwise bounds on $T^i\Omega^j S^k P_\sigma^{-1} g$ using Sobolev embeddings and Proposition \ref{Lemma:ijkresest}. The same result with different regularity assumptions can be found in \cite{Ta:13} (cf.\  Proposition 16).

Let $\Gamma$ denote the collection of all vector fields in $\Omega$,
$T$, and $S$. We write $\Gamma^\alpha$ to denote a product of these vector fields indexed by the multiindex $\alpha$ and let  $\Gamma^{\le n}$ denote a linear combination of
$\Gamma^\alpha$ for $|\alpha| \le n: \Gamma^{\le n} :=
\sum_{|\alpha| \le n}c_\alpha \Gamma^\alpha$. For the sake of notational
simplicity, we write $v_{ijk} := T^i \Omega^j S^k v$ and $g_{ijk} :=
T^i\Omega^jS^kg$. Similarly we write $v_{\le i \le j \le k} :=
T^{\le i}\Omega^{\le j}S^{\le k}v$ and  $g_{\le i \le j \le k} :=
T^{\le i}\Omega^{\le j}S^{\le k}g$. 

	\begin{lemma} \label{Lemma:vijkptbnd}
		Let $\Im \sigma \ge 0$ with $|\sigma| \gtrsim 1$ and assume $g \in \mathcal{LE}^*$ satisfies \eqref{gijkbnd}. 
		
		\begin{enumerate}
			\item[(i)] Then
		\begin{equation} \label{vijkptbnd}
			|T^i\Omega^jS^k P_\sigma^{-1}g| \lesssim (|\sigma|\langle r \rangle)^{-1}, \quad i + (\loss +1)j + (\loss +1)k \le M - 5(\loss+1)
		\end{equation}
		
	\item[(ii)] 	 If, in addition, we have $\sigma \in \mathbb{R}$ and $M \ge 4$, then we have the outgoing radiation condition 
			\begin{equation}\label{radcond2}
				\lim_{r \to \infty} r(\pa_r-i\sigma)T^i\Omega^jS^kv(\sigma) = 0, \quad i + (\loss +1)j + (\loss +1)k \le M - 5(\loss+1).
			\end{equation} 
		
		\end{enumerate}
	\end{lemma}
	
	\begin{proof}
		(i) To begin, we claim that establishing the estimate
        \begin{equation} \label{keyest}
		 		\sum_m 2^{\frac{m}{2}}\|(\partial_r^2 + \tau^2)( r v_{ijk})\|\lrsupo \lesssim 1 
			\end{equation}
	will suffice to obtain \eqref{vijkptbnd}.  To see this, note that the fundamental solution to  $(\partial_r^2+\tau^2)$ (for $\Im \sigma >0$) is given by $\phi_\sigma(s) = \sigma^{-1}e^{-i\sigma |s|}$. Thus (extending by continuity to $\Im \sigma \geq 0$), $rv_{ijk} = (\partial_r^2+\tau^2)(rv_{ijk}) * \phi_\sigma$.	Splitting
	$$
	\smallabs{r v_{ijk}}\leq \sum_m \int_{A_m} \smallabs{(\pa_r^2+\sigma^2)(r v_{ijk})(s)\phi_\sigma(r-s)} \, ds
	$$
	we now use the fact that $|\phi_\sigma(s)|\lesssim |\sigma|^{-1}$ and apply the Cauchy--Schwarz inequality to each term in the sum to 
	 to find
	$$
		|\sigma||rv_{ijk}| \lesssim |\sigma| \cdot |\sigma|^{-1} \sum_{m}2^{\frac{m}{2}}\|(\partial_r^2+\tau^2)(rv_{ijk})\|_{L^2_rL_\omega^\infty(A_m)} \lesssim 1.
	$$
         This yields \eqref{vijkptbnd} with $\la r \ra$
          replaced by $r,$
          i.e.,
          		\begin{equation} \label{vijkptbnd01}
			| rT^i\Omega^jS^k P_\sigma^{-1}g| \lesssim |\sigma|, \quad i + (\loss +1)j + (\loss +1)k \le M - 5(\loss+1).
		\end{equation}

          That the estimate in fact holds for
          $\la r \ra $ essentially follows from the fact that there is
          no preferred origin to our coordinate system, and our
          estimates are translation-invariant.  In particular, letting
          $\realtau$ denote the translation operation $\realtau f(x)=f(x-a)$
          with $a$ a fixed nonzero vector, we observe that $\realtau$
          commutes with translation vector fields $T_i,$ while
          $$
          \realtau\Omega-\Omega \realtau,\ \realtau S-S \realtau
          $$
          are both of the form $\sum c_j T_j \realtau=\sum c_j \realtau T_j.$
          Thus
          $$
T^i \Omega^j S^k \realtau-\realtau T^i \Omega^j S^k =\sum \Gamma^\alpha \realtau
=\sum \realtau \widetilde{\Gamma}^\alpha
$$
where sums are over products of vector fields $\Gamma^\alpha,\
\widetilde{\Gamma}^\alpha$ of the form $T^{i'} \Omega^{j'} S^{k'}$
with $i'+(\loss+1)j' + (\loss+1)k'\leq i+(\loss+1)j + (\loss+1)k.$

Note further that the hypotheses on $P$ are translation-invariant, so
that if $P_\sigma v=g$ then $P'_\sigma \realtau v=\realtau g$ with $P'_\sigma$
an operator satisfying the same hypotheses.  Owing to the commutation
properties of $\realtau$ with our rest operators, $\realtau g$ also satisfies \eqref{gijkbnd}, hence $\realtau
v$ likewise satisfies \eqref{vijkptbnd01}.  Now translating back and again using the commutation
properties of the vector fields $\Gamma^\alpha$ shows that
		\begin{equation} \label{ick}
			|(x-a) T^i\Omega^jS^k v| \lesssim |\sigma|, \quad i + (\loss +1)j + (\loss +1)k \le M - 5(\loss+1).
                      \end{equation}
Adding \eqref{vijkptbnd01} and \eqref{ick} yields the
                    desired estimate.

	To establish \eqref{keyest} we write
	$$
				(\partial_r^2+\sigma^2) = -P_\sigma -(2r^{-1}\partial_r + r^{-2} \Delta_\angvar + Q)
			$$
		and commute $P_\sigma$ with $T^i\Omega^jS^k$ to
                find
		\begin{equation} \label{drtaubnd}
			(\partial_r^2+\sigma^2)(rv_{ijk}) = -r^{-1} \Delta_\angvar v_{ijk} + rQv_{\le i\le j\le k} - r g_{\le i \le j \le k}
		\end{equation}
		and bound each term on the right hand side as in \eqref{keyest}.
	
	To see where the vector field loss occurs, consider the term $r^{-1}\Delta_\angvar v_{ijk}$ on the right hand side of \eqref{drtaubnd}. We wish to show this term satisfies \eqref{keyest}. By the spherical Sobolev embedding, it suffices to show
	\begin{equation} \label{goal} 
		\sum_m  \| \la r \ra^{-\frac{3}{2}} \Delta_\angvar v_{ijk} \|\ltwoam + \| \la r \ra^{-\frac{3}{2}}  \Omega^2(\Delta_\angvar v_{ijk}) \|\ltwoam \lesssim 1 .
	\end{equation}
	We have 
		\begin{equation} \label{largeletau}
			\| \la r \ra^{-\frac{1}{2}} v_{ijk} \|_{L^2(A_m)} \lesssim |\tau| \| \la r \ra^{-\frac{1}{2}} v_{ijk} \|\ltwoam \lesssim \|v_{ijk}\|\letn{} \lesssim 1
		\end{equation}
when $ i + (\loss+1)j + (\loss+1)k < M-\loss-1$. 

 	Replacing $\Delta_\angvar$ by $\sum_{|\alpha| = 2} \Omega^\alpha$ and using \eqref{largeletau} then yields
 	\begin{align*} 
 		LHS \hbox{ of } \eqref{goal} &= \sum_m \| \la r \ra^{-\frac{3}{2}} v_{i(j+2)k} \|\ltwoam + \| \la r \ra^{-\frac{3}{2}} v_{i(j+4)k} \|\ltwoam \\
 			&\lesssim 1 
 	\end{align*}	
when $i + (\loss+1)j + (\loss+1)k \le M - 5(\loss+1)$, as desired. The remaining terms on the right hand side of \eqref{drtaubnd} are handled similarly.

 (ii) Note that we have
	\[ (\partial_r + i\tau)(\partial_r - i\tau)(rv_{ijk}) = (\partial_r^2 + \tau^2)(rv_{ijk}) \]
	so that
	\begin{equation} 
		\partial_r\big(e^{i\tau r} (\partial_r - i\tau)(rv_{ijk})\big)= (\partial_r^2 + \tau^2)(rv_{ijk})e^{i\tau r}. \label{odesol}
	\end{equation}
By \eqref{keyest} we see that $(\pa_r^2 + \sigma^2)(rv_{ijk})$ is integrable in $r$, so that the limit $\lim_{|x| \to \infty} |(\partial_r - i\tau)(rv_{ijk})|$ exists for each $\theta \in \mathbb{S}^2$  since
	\[ \lim_{|x| \to \infty} |(\partial_r -i\tau)(rv_{ijk})| = \Big| \int_0^\infty (\partial_r^2 + \tau^2)(rv_{ijk}) e^{i\tau r} \dd r + v_{ijk}(0)\Big|. \]
For fixed $\theta \in \mathbb{S}^2$, take
	\[ c_\theta = \lim_{|x| \to \infty} |(\partial_r -i\tau)(rv_{ijk})| = \lim_{|x| \to \infty} |r(\partial_r - i\tau)v_{ijk} + v_{ijk}|. \]
	By part (i) of this proposition, $\lim_{r \to \infty} v_{ijk} = 0$ so that 
		$$\lim_{r \to \infty} |r(\partial_r - i\tau)v_{ijk}| = c_\theta.$$ 
		Thus we can write $|(\partial_r - i\tau)v_{ijk}| = c_\theta r^{-1} + o(r^{-1})$.
	
	On the other hand, since $\tau$ is real and $5(\loss +1) \ge 4$, by Lemma \ref{Lemma:radcond} $v_{ijk}$ satisfies the radiation condition \eqref{radcond}. Then Sobolev embedding  yields
	\begin{equation} \label{polarradcond}
		\lim_{m \to \infty} 2^{\frac{m}{2}} \| (\partial_r - i\tau) v_{ijk} \|_{L_r^2L_\theta^\infty(A_m)} = 0.
	\end{equation} 
	It follows that 
	\[ 0  \ge \lim_{m \to \infty} 2^{\frac{m}{2}} \| c_\theta r^{-1} + o(r^{-1}) \|_{L^2_r(A_m)} = \lim_{m \to \infty} 2^{\frac{m}{2}} c_\theta 2^{-\frac{m}{2}} = c_\theta  \]
	Thus $c_\theta \equiv 0$ for $\theta \in \mathbb{S}^2$ so that $(\partial_r - i\tau)v_{ijk} \in o(r^{-1})$, which concludes the proof of \eqref{radcond2}.

	\end{proof}

       Finally we provide the following result on the finite order conormal regularity of $P_\sigma^{-1}$ for large $\sigma$ (see \cite{Ta:13} Proposition 17 and \cite{Mo:20} Proposition 6.5). We give a sketch of the proof, with a focus on the numerology in the proposition.
       
	\begin{proposition} \label{largesigmaresbnd}
		Let $|\sigma| \gtrsim 1$ with $\sigma$ real and take $g \in H_b^{s,\ell}$ with $\ell > \frac{1}{2}$. Then
		$$
			\left| (\sigma \partial_\sigma)^p (e^{-i\sigma r}P_\sigma^{-1}g ) \right| \lesssim |\sigma|^{p-1}\langle r \rangle^{-1+p(1-\varepsilon)}, \quad p \le \ell - \frac{1}{2}, \quad s \ge (2p + 5)(\loss + 1)
		$$
for some $\epsilon \in (0,1]$.
	\end{proposition}
	
	\begin{proof}
		Identifying the vector fields $T, \Omega, r\partial_r$
                with b-vector fields and using the change of
                coordinates $\rho = r^{-1}$, we see $g \in
                \Hb^{s,\ell}$ implies
                $$\|r^{\ell-\frac{1}{2}-}T^i\Omega^jS^k g\|\les
                \lesssim1$$ for $i+j+k\leq s$. 

Take $v = P_\sigma^{-1} g$ and write
		\begin{equation}		
		\begin{split} 
			(\tau \partial_\tau)^p\left(ve^{-i r  \tau}\right) &= \big( (-S+r(\partial_r-i\sigma))^p v \big)e^{-ir\sigma}\\
			&= \left(\sum_{m=0}^p c_m [r(\partial_r-i\sigma)]^m(-S)^{p-m}v \right)e^{-ir\sigma}\\
			&= \left(-S^pv + \sum_{m = 1}^p \sum_{n = 1}^m c_{mn} r^n(\partial_r-i\tau)^n(-S)^{p-m}v \right)e^{-i  r  \tau}. 
		\end{split}
		\end{equation}
We claim $|( \partial_r - i\tau)^p v_{ijk}| \lesssim |\sigma|^{p-1}r^{-1-p\varepsilon} $ with $\varepsilon \in (0,1]$ for 
	$$
		i + (\loss +1)j + (\loss+1)k \le s-5(\loss +1) - 2(\loss +1)p.
	$$ 
The proposition follows from the claim. To see this, assume the claim holds. Then $|r^n (\partial_r-i\sigma)^n(-S)^{p-m}v| \lesssim |\sigma|^{n-1}r^{-1+n(1-\varepsilon)}$ when 
	$$
		(\loss+1)p \le s-5(\loss +1) + (\loss+1)(m-2n)
	$$ 
for $0 \le n \le m \le p$. The smallest value of $ s-5(\loss +1) + (\loss+1)(m-2n)$ in this range is $s-(5+p)(\loss+1)$, which justifies the assumption $s \ge (2p+5)(\loss + 1)$.

	It is left to prove the claim. When $p=0$, the claim follows from \eqref{vijkptbnd}. For general $p$, applying $(\partial_r-i\tau)^{p-1}$ to \eqref{drtaubnd} we find 
	\begin{equation} \label{gencase}
	\begin{split} 
	(\partial_r &+ i\tau)(\partial_r - i\tau)^p(rv_{ijk})\\
		&= \sum_{m=0}^{p-1}\left((-1)^{p-m+1}c_m r^{-(p-m)} (\partial_r - i\tau)^m \Delta_\angvar  v_{ijk}\right) + r(\partial_r - i\tau)^{p-1} Q v_{\le i \le j \le k}  \\
			&\quad + C(\partial_r - i\tau)^{p-2} Q v_{\le i \le j \le k} + \big(r(\partial_r - i\tau)^{p-1} + C(\partial_r - i\tau)^{p-2}\big) g_{\le i\le j\le k}.  
	\end{split}	
	\end{equation}
Each $\vijk$ term on the right hand side of \eqref{gencase} is
  bounded in magnitude by $|\tau|^{p-1} \la r \ra^{-1-p\varepsilon}$ inductively
  using \eqref{vijkptbnd}. We note the $\varepsilon$ factor in the exponent shows up when considering, for example, the $\sigma S^{-\kappa}\partial_x$ term in $Q$ which yields a $S^{-\kappa+1}(\partial_r-i\sigma)^{p-1}v_{i+1,j,k}$ term on the right hand side of \eqref{gencase}\footnote{When $\kappa>2$ we have $\varepsilon =1$ since in this case $S^{-\kappa+1} \subseteq S^{-1}$. When $\kappa \in (1,2)$, it should be possible to improve the estimate to $\varepsilon=1$. We roughly illustrate the idea of the process by considering $$(\sigma \partial_\sigma)(ve^{-ir\sigma}) = \big([-S+r(\partial_r-i\sigma)]v\big)e^{-i\sigma r} $$ in which case we have
	$$
		(\partial_r^2 + \sigma^2)(rv) = -r^{-1} \Delta_\angvar v + rQv -rg
	$$
and the $S^{-\kappa}$ coefficient terms in $Q$ do not satisfy the $r^{-2}$ pointwise bounds. In this case we write
	$$
		(\partial_r+i\sigma)(\partial_r-i\sigma)(rv) = S^{-\kappa}(\partial_r-i\sigma)(rv) + \mathcal{O}(r^{-2})
	$$
then solve for $(\partial_r-i\sigma)(rv)$ by integrating as above to get $(\partial_r-i\sigma)v \lesssim r^{-2}$, as needed.}. This requires 
	$$
		i + (\loss +1)j + (\loss+1)k \le s-5(\loss +1) - 2(\loss +1)p.
	$$ 
	The loss of $2(\loss +1)p$ follows from the $\Delta_\angvar$
        operator on the right hand side of \eqref{gencase}, which we
        replace by $\Omega^2$. 
	
	To handle the $g_{\le i\le j\le k}$ terms we note $g_{\le i\le j\le k} \in \Hb^{s-i-j-k,\ell}$ and by Sobolev embeddings we have $|\partial_r^p g_{\le i\le j\le k}| \lesssim \langle r \rangle^{-\ell-\frac{3}{2}-p}$ if $s-i-j-k-p>3$, which indeed holds from the above restrictions on $i,j,k$. Now we calculate
	\begin{align*}
		|r&(\partial_r-i\tau)^{p-1}g_{\le i \le j \le k}| + |(\partial_r-i\tau)^{p-2}g_{\le i \le j \le k}|\\ 
			&= \big|r \sum_{m = 0}^{p-1} c_m (-i\tau)^{m}\partial_r^{p-1-m}g_{\le i \le j \le k}\big| + \big|\sum_{m = 0}^{p-2}(-i\tau)^{m}c_m \partial_r^{p-2-m}g_{\le i \le j \le k}\big| \\
			&\lesssim |\sigma|^{p-1} \la r \ra^{-\ell-\frac{1}{2}}
	\end{align*}
	which yields the desired $|\tau|^{p-1} \la r \ra^{-1-p}$ bound when $p \le \ell - \frac{1}{2}$.
	
We write
	$$ 
		\partial_r[((\partial_r-i\tau)^prv_{ijk})e^{i\tau r}] = [(\partial_r+i\tau)(\partial_r-i\tau)^p(rv_{ijk})]e^{i\tau r}.
	$$
	When $p = 1$, Lemma \ref{Lemma:vijkptbnd}  shows $\lim_{r \to \infty} (\pa_r-i\sigma)(rv_{ijk}) =0$, so we can integrate from infinity to prove the claim.   When $p >1$ note that 
	\[ (\partial_r - i\tau) v_{ijk} =   r ^{-1} x\cdot \nabla v_{ijk} - i\tau v\]
	and thus
	\begin{align*} 
		|(\partial_r - i\tau)^2 v_{ijk}| &= |x  r ^{-1} (\partial_r - i\tau) v_{(i+1)jk} - i\tau (\partial_r - i\tau) v| \\
			& \lesssim |(\partial_r -i\tau) v_{(i+1)jk}| + |\tau||(\partial_r-i\tau)v|. 
	\end{align*}
	Iterating, we find 
	$$
	 |(\pa_r-i\sigma)^p v_{ijk}| \lesssim \sum_{m=0}^{p-1} |\sigma|^m|(\pa_r-i\sigma)v_{(i+p-1-m)jk}|
	$$
	and thus $r(\pa_r-i\sigma)^pv_{ijk} \to 0$ as $r \to \infty$ by \eqref{radcond2} for
	$$ i +(\loss+1)j +(\loss+1)k \le s -5(\loss +1) - (p-1) $$
	which is satisfied under the conditions of the claim since $-2(\loss+1)p \le -p+1$. Thus when $p>1$ we can integrate from infinity just as in the $p=1$ case to prove the claim.

	\end{proof}

\section{Low Frequency Conormal Estimates} \label{LowEnergy}
The purpose of this section is to provide information about the behavior of $(\sigma\partial_\sigma)^MP_\sigma^{-1}f$ for $f$ in an appropriately chosen function space.  We follow the approach introduced by Hintz in \cite{Hi:20}.

In order to obtain good asymptotic expansions as $\sigma \to 0,$ it is useful to work with the conjugated operator
	$$
		\tP(\sigma) \equiv e^{-ir\sigma}P_\sigma e^{ir\sigma}.
	$$
For $\sigma \in \RR,$ we define $\tP(\sigma)^{-1}$ to be $\lim_{\ep \downarrow 0} \tP(\sigma+i\varepsilon)^{-1}.$ Note that $\tP(0) = P_0$.

Working on the spatial manifold $X$, we find
\[
\begin{split}
	\tP(\sigma) &= \rho^2[-(\rho\partial_\rho)^2+\rho\partial_\rho+\Delta_\angvar] + 2i\rho\sigma(\rho\partial_\rho-1)+\A^\kappa(X)\sigma^2 \\
		& \quad +\A^{\kappa+1}(X)\sigma \Diffb^1(X)+\A^{\kappa+2} \Diffb^2(X),
\end{split}
\]
where the operator class
$$
\A^{\kappa+2} \Diff_b^2(X)
$$
refers to b-differential operators with conormal coefficients of the
specified order.

Note that
\begin{equation}\label{tPdef}
\tP(\sigma) = \tP(0)-\sigma R
\end{equation}
where
\begin{equation}\label{Rdef}
R= -2i\rho(\rho\partial_\rho-1)+\A^{\kappa+1}\Diff_b^1(X) + \sigma \A^\kappa.
\end{equation}
In Proposition \ref{prop:lowenergy} we study the low frequency behavior of $\tP(\sigma)^{-1} f$ by using a formal Neumann series argument to write
$$
	\tP(\sigma)^{-1} = \tP(0)^{-1}\big( \Id +\dots + (\sigma R \tP(0)^{-1} )^N \big) +
    \tP(\sigma)^{-1} (\sigma R \tP(0)^{-1} )^{N+1}.
$$
The key feature of this expansion is that each iterative application of $R\tP(0)^{-1}$ results in the loss of one power of decay. After enough iterations, the output will be too large to apply $\tP(0)^{-1}$ again, which forces the iteration to stop with the final application of $\tP(\sigma)^{-1}$.

In subsection \ref{subsection:Vasy} we record a previously established preliminary estimate which will aid our subsequent calculations. We then prove the necessary mapping properties for $\tP(0)^{-1}$ and $\tP(\sigma)^{-1}$ in subsection \ref{subsection:mapprop}. This includes an expansion for $\tP(0)^{-1}f$ in Lemma \ref{lemma:P0inverse} (and in Lemma \ref{lemma:P0inverselessdecay} for larger inputs). We note that in Lemma \ref{lemma:P0inverse}, we see that even if $f$ is rapidly decaying, the expansion for $P(0)^{-1}$ is limited by the perturbative $\rho^\kappa L_1$ term in $\tP(0)$. It is this limitation which dictates when the Neumann series must end, which in turn dictates the final decay rate obtained in the proof of the main theorem. Finally, in subsection \ref{subsection:conormalestimates} we study the Neumann series and obtain the desired conormal estimates.

\subsection{Preliminary results} \label{subsection:Vasy}
The following result from \cite{Va:20} (cf.\ Theorem 1.1) will be used to analyze the low energy resolvent:
\begin{theorem}\label{theorem:lowerbnd}
	For $s,\ell,\nu \in \mathbb{R}$ with $\ell<-\frac{1}{2}, s+\ell>-\frac{1}{2}, \ell-\nu \in (-\frac{3}{2},-\frac{1}{2})$ the bound
	$$
		\|(\rho+|\sigma|)^\nu u\|_{H_b^{s,\ell}(X)} \lesssim \|(\rho + |\sigma|)^{\nu-1}\tP(\sigma)u\|_{H_b^{s,\ell+1}(X)}
	$$
holds for bounded $\sigma$.
\end{theorem}
Note that this result is obtained from taking $r=s+l$ in Theorem~1.1
of \cite{Va:20} to obtain
$$
\smallnorm{(\rho+\smallabs{\sigma})^\nu
  u}_{H_{\scc,b,\res}^{s,r,l}}\lesssim
\smallnorm{(\rho+\smallabs{\sigma})^{\nu-1}
  \tP(\sigma) u}_{H_{\scc,b,\res}^{s-2,r+1,l+1}};
$$
these triple-index Sobolev spaces measure a combined (``second
microlocal'') b and scattering regularity.
We then coarsen the estimate by estimating the right hand-side by
$$
\smallnorm{(\rho+\smallabs{\sigma})^{\nu-1}
  \tP(\sigma)  u}_{H_{\scc,b,\res}^{s,r+1,l+1}};
$$
Since $H_{\scc, b,\res}^{s,s+l,l} = H_b^{s,l}$ the result then follows.

We remark here that giving up two derivatives as we do in
Theorem~\ref{theorem:lowerbnd} in return for weighted estimates does
have consequences for the numerology of our decay hypotheses on
initial data (mainly owing to the conjugating factor
$e^{i\sigma/\rho}$ which means that non-Schwartz data has limited
b-regularity after multiplication by this factor).
It is possible that a finer accounting of
regularity in our iteration, tracked in the second microlocal Sobolev
spaces of \cite{Va:20}, would yield more precise decay hypotheses.

We use this theorem to establish basic mapping properties of $\tP(0)^{-1}$ that will help us find a useful expression for $\tP(0)^{-1}f$.

\begin{corollary} \label{cor:prelimmapprop}
Let $\ell \in (-\frac{3}{2}, -\frac{1}{2})$ and $s+\ell > -\frac{1}{2}$. Then
	$$
		\tP(0)^{-1}: H_b^{s,\ell+2} \to H_b^{s,\ell}.
	$$
\end{corollary}

\subsection{Low frequency mapping properties}\label{subsection:mapprop}
Following Hintz \cite{Hi:20}, we write
$$
\tP(0)= \rho^2(L_0 + \rho^\kappa L_1)
$$
with
$$
L_0=-(\rho\pa_\rho)^2 +\rho \pa_\rho + \Lap_\angvar,
$$
$$
L_1 \in \A^0 \Diffb^2.
$$

The conclusion of the following lemma involves sums over
finite-dimensional spaces of spherical harmonics, schematically
denoted
$$
\sum_{j=1}^M \rho^j Y_{j-1}
$$
Since all Sobolev norms are equivalent on such spaces, we write
$$
\abs{Y_{j-1}}
$$
to denote the supremum of each component, but will use
tacitly the fact that this is equivalent to taking any desired
b-Sobolev norm of these angular pieces, as well.

We begin by considering $L_0^{-1}$ then argue perturbatively to analyze $\tP(0)^{-1}$.
\begin{lemma}\label{lemma:L0inverse}
	Let $f \in H_b^{s,\gamma}$ with $\gamma > -\frac{1}{2}$ and $s>0$. Assume $u \in H_b^{s,-\frac{1}{2}-}$ solves $L_0 u = f$. If $\gamma + \frac{3}{2} \notin \mathbb{N}$ then
	$$
		u = \sum_{j=1}^{\floor{\gamma+\frac{3}{2}}} \rho^j Y_{j-1} + q
	$$
	where $Y_j$ is a linear combination of $j^{th}$ order
        spherical harmonics, $q \in H_b^{s+2,\gamma-},$ and
	$$
		\sum_{j=1}^{\floor{\gamma+\frac{3}{2}}} |Y_{j-1}| + \|q\|_{H_b^{s+2,\gamma-}} \lesssim \|f\|_{H_b^{s,\gamma}}.
	$$
\end{lemma}

\begin{proof}
	We are interested in the behavior of $u$ as $\rho \to 0$.  Take $\chi_\partial(\rho)$ to be a cutoff which is 1 on a neighborhood of $\partial X$ and 0 for $\rho \ge 1$. Define $u_\partial := \chi_\partial u$. Then we have
	\begin{equation}\label{L0ubndry}
	L_0 u_\partial = f_\partial \in H_b^{s,\gamma}(X)
	\end{equation}
with $\hbox{supp}(f_\partial) \subset [0,1)_\rho$.

	The Mellin transform in $\rho$ of a function $g$ defined on
        $X$ is given by $\M {g}(\xi) := \int_0^\infty \rho^{-i\xi}
        g(\rho) \frac{d\rho}{\rho}$. Taking the Mellin transform in
        $\rho$ of Equation \eqref{L0ubndry} yields 
	$$
		\hat{L}_0(\xi) \M {u}_\partial(\xi)=\M {f}_\partial(\xi), \quad \hat{L}_0(\xi) = -(i\xi)^2+ i\xi + \Delta_\angvar.
	$$
        The Mellin transform of $f_\partial$ is
          holomorphic in $\Im\xi > - \gamma - \frac{3}{2}$
          (see Proposition \ref{prop:mellin} in Appendix
          \ref{app:bgeometry}). Since $u \in H_b^{s,-\frac{1}{2}- }$
          by assumption, we can invert the Mellin transform of
          $u_\partial$ by integrating along a contour in Im$\xi >
          -1$. Thus we obtain
	\begin{equation} \label{invmel}
		u_\partial(\rho) = \int_{\Im{\xi} = -1+\varepsilon}\rho^{i\xi} \M {u}_\partial(\xi) \, d\xi.
	\end{equation}
	
For each $\xi \in \mathbb{C}$ we decompose $\M {u}_\partial(\xi)$ and $\M {f}_\partial(\xi)$ into spherical harmonics: $\M {u}_\partial(\xi) = \sum_{j=0}^\infty \sum_{m=-j}^j \M {u}_{mj}(\xi) y_{mj}$ where $$\M {u}_{mj}(\xi) = \int_{\mathbb{S}^2} \M {u}(\xi)y_{mj} \, d\angvar$$
and similarly for $\M {f}_\partial$. Define $\mathcal{Y}_j := \hbox{span}_m(y_{mj})$ and note $\hat{L}_0(\xi)\big|_{\mathcal{Y}_j} = -(i\xi)^2+i\xi+j(j+1)$. Thus $\hat{L}_0(\xi)^{-1}\big|_{\mathcal{Y}_j} = \frac{1}{-(i\xi)^2+i\xi+j(j+1)}$ has simple poles at $\xi = ij, - i(j+1)$. Using the spherical harmonic decomposition we have
	$$
		\M {u}_\partial(\xi) = \sum_{j=0}^\infty \sum_{m=-j}^j \hat{L}_0(\xi)^{-1}\big|_{\mathcal{Y}_j} \M {f}_{mj}(\xi) y_{mj}
	$$
Since $f_\partial \in H_b^{s,\gamma-}$, we can push the contour of
integration in \eqref{invmel} down to $\Im{\xi} =
-\gamma-\frac{3}{2}+\varepsilon$ and pick up residues at $\xi = -i,
-2i, \dots, -\floor{\gamma+ \frac{3}{2}}i$. The residue of
$\rho^{i\xi} \M {u}_{k-1}(\xi)$  at $\xi = -ik$ is $$\rho^k
\sum_{m=-k+1}^{k-1} i(2k-1)^{-1} \M {f}_{m(k-1)}(-ik)y_{m(k-1)}.$$ Thus
	$$
		u_\partial(\rho,\angvar) = \sum_{k=1}^{\floor{\gamma+
                    \frac{3}{2}}} \rho^k Y_{k-1} + \int_{\Im{\xi} =
                  -\gamma-\frac{3}{2}+\varepsilon} \rho^{i\xi} \M
                {u}_\partial(\xi) d\xi 
	$$
where $Y_{k-1} = \sum_{m=-k+1}^{k-1} i(2k-1)^{-1} \M
{f}_{m(k-1)}(-ik)y_{m(k-1)}$ is a linear combination of $(k-1)^{th}$
order spherical harmonics. By Cauchy-Schwarz, $$|\M
{f}_\partial(-ik)|^2 \lesssim \int_0^\infty \rho^{-2\gamma}
|f_\partial|^2 \frac{d\rho}{\rho^{-4}}$$ when $k < \gamma +
\frac{3}{2}$. Thus  
	$$
		|\M {f}_{m(k-1)}(-ik)|^2 \le \int_{\mathbb{S}^2} |\M {f}_\partial(-ik) y_{m(k-1)}| \, d\angvar \lesssim \|f_\partial \|_{H_b^{s,\gamma}}
	$$
and $|Y_{k-1}| \lesssim \|f_\partial \|_{H_b^{s,\gamma}}$ as
desired. 
 Finally we have obtained $$q = \int_{\Im{\xi} =
  -\gamma-\frac{3}{2}+\varepsilon} \rho^{i\xi} \hat{L}_0(\xi)^{-1} \M
{f}_\partial(\xi) d\xi.$$  That $\M q$ is holomorphic in
$\Im (\bullet) >-\gamma-3/2+\ep$ follows since $q$ differs from
$$
\M^{-1}(\hat{L}_0(\xi)^{-1} \M f_\pa(\xi))
$$
by subtraction of the poles of $\hat{L}_0(\xi)^{-1} \M f_\pa(\xi)$ in
this region.  Consequently, since
$\hat{L}_0(\xi)\rvert_{\mathcal{Y}_j}$ grows quadratically in both
$\xi$ and $j,$ for any $\mu>-\gamma-3/2+\ep$ and for $\nu
\in \RR,$
$$
\abs{\M q_j(\nu+i\mu)}\leq C (1+j^2+\nu^2)^{-1} \abs{\M f_{j}(\nu+i\mu)},
$$
hence
$$
\norm{(\xi^2+\Lap_\theta)^{s+2}\M q(\nu+i\mu)}_{L^2_\nu}^2\leq C \norm{(1+\xi^2+\Lap_\theta)^{s/2}\M f_{\pa}(\nu+i\mu)}_{L^2_\nu}^2
$$
and the estimate $\|q\|_{H_b^{s+2,\gamma-}} \lesssim
\|f_\partial\|_{H_b^{s,\gamma}}$ follows, using Proposition \ref{prop:mellin}.
\end{proof}

The following lemma is analogous to Lemma \ref{lemma:L0inverse} but with less decay assumed for $f$. It is stated separately because of minor technical changes in the numerology (note we now assume $u \in H_b^{s,\gamma-}$). The proof is the same and we provide an abbreviated argument. 

\begin{lemma}\label{lemma:L0inverselessdecay}
	Let $f \in H_b^{s,\gamma}$ with $\gamma \in (-\frac{3}{2},-\frac{1}{2})$. Assume $u \in H_b^{s,\gamma-}$ solves $L_0 u = f$. Then
	$$
		\|u\|_{H_b^{s+2,\gamma-}} \lesssim \|f\|_{H_b^{s,\gamma}}.
	$$
\end{lemma}

\begin{proof}
	The lemma is proved as in Lemma \ref{lemma:L0inverse}. Instead of \eqref{invmel} the assumptions on $u$ allow us to integrate along $\Im \xi = -\gamma-\frac{3}{2}+\epsilon$. Note $-\gamma-\frac{3}{2}+\epsilon \in (0,-1)$. Since $f$ is assumed to have the same decay, we cannot push the contour of integration further into the lower half plane. The result follows immediately using Proposition \ref{prop:mellin} as before.
\end{proof}

The preceding lemmas will now be used to establish asymptotic expansions for $\tP(0)^{-1}f$ which depend on the amount of decay assumed for $f$.

\begin{lemma}\label{lemma:P0inverse}
	Let $f \in H_b^{s,\gamma}$ for $s>0$ and $\frac{3}{2} < \gamma \le \frac{3}{2} + \kappa$ with $\gamma + \frac{3}{2} \notin \mathbb{N}$. Then
	\[ \tP(0)^{-1}f = \sum_{j=1}^{\floor{\gamma - \frac{1}{2}}}
          \rho^{j}Y_{j-1} + q \]
        with
        $$
q\in H_b^{s+2,\gamma-2-}
        $$
	and
	\[\sum_{j=1}^{\floor{\gamma - \frac{1}{2}}} |Y_{j-1}| + \|q\|_{H_b^{s+2,\gamma-2-}} \lesssim \|f\|_{H_b^{s,\gamma-}}. \]
\end{lemma}

\begin{proof}
	Define $u := \tP(0)^{-1}f$ and note that $u \in
        H_b^{s,-\frac{1}{2}-}$ by Corollary
        \ref{cor:prelimmapprop}. Writing $L_0 = \rho^{-2} \tP(0) -
        \rho^\kappa L_1$ with $L_1 \in \A^0 \Diff_b^2$, we find $L_0u =
        \rho^{-2}f + \rho^\kappa L_1 u$. Then by Lemma
        \ref{lemma:L0inverse} we have   
	\[u = \sum_{j=1}^{\floor{\gamma-\frac{1}{2}}} \rho^j Y_{j-1} + q\]
with
	\begin{align*} 
		\sum_{j=1}^{\floor{\gamma-\frac{1}{2}}} |Y_{j-1}| + \|q\|_{H_b^{s+2,\gamma-2-}} & \lesssim \left\| \rho^{-2}f + \rho^\kappa L_1\left(\sum_{j=1}^{\floor{\gamma-\frac{1}{2}}} \rho^j Y_{j-1} + q\right)\right\|_{H_b^{s,\gamma-2-}}\\
			&\lesssim \|f\|_{H_b^{s,\gamma-}} + \varepsilon\Big(\sum_{j=1}^{\floor{\gamma-\frac{1}{2}}} |Y_{j-1}| + \|q\|_{H_b^{s+2,\gamma-2-}}\Big)
	\end{align*}
	for $\rho$ small for some $0< \varepsilon < 1$.  (Recall that
        all norms are comparable to one another on the finitely many
        angular modes $Y_j$.) Bootstrapping the last terms on the right hand side above to the left hand side then yields the desired inequality.
\end{proof}

Note even if $f$ has faster decay as $\rho \to 0$ than assumed in Lemma \ref{lemma:P0inverse}, there is no improvement over the result for $f \in H_b^{s,\kappa+\frac{3}{2}}$ due to the perturbative $\rho^\kappa L_1$ term in $\tP(0)$. 

As before, we prove a lemma analogous to Lemma \ref{lemma:P0inverse} now assuming $f$ has less decay. More regularity is also assumed for $f$ due to the numerology in Corollary \ref{cor:prelimmapprop}.

\begin{lemma}\label{lemma:P0inverselessdecay}
	Let $f \in H_b^{s,\gamma}$ for $s>1$ and $\gamma \in (\frac{1}{2},\frac{3}{2})$ with $\gamma + \frac{3}{2} \notin \N$.  Then $u = \tP(0)^{-1}f$ satisfies
	$$
		\|u\|_{H_b^{s+2,\gamma-2-}} \lesssim \|f\|_{H_b^{s,\gamma}}.
	$$
\end{lemma}

\begin{proof}
	The proof is analogous to that of Lemma \ref{lemma:P0inverse} except Corollary~\ref{cor:prelimmapprop}  now implies $u \in H_b^{s,\gamma-2-}$. Note $\gamma - 2 \in (-\frac{3}{2},-\frac{1}{2})$. As before we find $L_0u = \rho^{-2}f - \rho^\kappa L_1u$. Now we use Lemma \ref{lemma:L0inverselessdecay} to find
	$$
		\|u\|_{H_b^{s+2,\gamma-2-}} \lesssim \| \rho^{-2}f\|_{H_b^{s,\gamma-2-}} + \|\rho^\kappa L_1u\|_{H_b^{s,\gamma-2-}} \lesssim \|f\|_{H_b^{s,\gamma-}} + \epsilon\|u\|_{H_b^{s+2,\gamma-2-}}
	$$
for $\rho$ small for some $0 < \epsilon <1$.
\end{proof}

Since $P(0)^{-1}f$ generates terms of the form $\rho^n Y_\ell$  and we will apply $P(0)^{-1}$ iteratively, we now consider the output of $P(0)^{-1} \rho^nY_\ell$.

\begin{lemma}\label{lemma:P0inversesh}
	If $3 \le n$ and $n \neq \ell+3$, then
	$$
	\tP(0)^{-1}(\rho^nY_\ell)= \rho^{n-2}Y_\ell + \sum_{j=1}^{\floor{\kappa}+1} \rho^jY_{j-1} + \A^{\kappa+1-}
	$$
\end{lemma}

\begin{proof}
	Note for general $m$ we have $L_0 \rho^mY_\ell = -(m+\ell)(m-\ell-1)\rho^mY_\ell$ so $L_0^{-1}(\rho^mY_\ell) = c\rho^mY_\ell$ when $m \neq -\ell,\ell+1$. 
	
	Define $u:=\tP(0)^{-1}(\rho^nY_\ell)$. By Sobolev embeddings \eqref{sobemb} we have $\rho^nY_\ell \in \A^n \subseteq H_b^{\infty,n-\frac{3}{2}-}$. It follows from Corollary \ref{cor:prelimmapprop} that $u \in H_b^{\infty,-\frac{1}{2}-}$.
	As before we use $L_0 = \rho^{-2}\tP(0)-\rho^\kappa L_1$ to write $L_0 u = \rho^{n-2}Y_\ell + H_b^{\infty,\kappa-\frac{1}{2}-}$. The result then follows by Lemma \ref{lemma:L0inverse} and using our assumptions on  $n$ to find $L_0^{-1}(\rho^{n-2}Y_\ell) = c\rho^{n-2}Y_\ell$.

\end{proof}

We now consider the mapping properties of $P(\sigma)^{-1}$, which will be needed in the last term in the Neumann series. The following lemma is implicit in the proof of Lemma 2.16 in
\cite{Hi:20}, and our proof follows that in \cite{Hi:20}.
\begin{lemma} \label{lemma:Psiginvlargeinput}
  Let $\alpha \in (0,1)$ and $s>1.$  Then there exists $\sigma_0$ such that $|\sigma| \le \sigma_0$ implies that for all $\delta>0$  sufficiently small,
  $$
\tP(\sigma)^{-1}\colon H_b^{s,1/2-\alpha} \to |\sigma|^{-\alpha-\delta} H_b^{s,-3/2-\delta}.
  $$
  \end{lemma}
  In \cite{Hi:20} this is used to show (by taking $s\to\infty$) that
  $$
\tP(\sigma)^{-1}\colon \A^{2-\alpha} \to |\sigma|^{-\alpha-0} \A^{-0}
$$
(and indeed, a more refined statement holds on the resolved space).
\begin{proof}
We take $l=-3/2-\delta$ in Theorem~\ref{theorem:lowerbnd}.  Then the
constraint on $s$ is $s>1+\delta,$ and is satisfied if $\delta>0$ is
sufficiently small.  The constraint on $\nu$ is $\nu \in
(-1-\delta,-\delta);$ we take $\nu=-2\delta$ to obtain
$$
\norm{(\rho+\smallabs{\sigma})^{-2\delta} u}_{s,-3/2-\delta}\lesssim
\norm{(\rho+\smallabs{\sigma})^{-1-2\delta} \tP(\sigma) u}_{s,-1/2-\delta}.
$$
Estimating
$$
(\rho+\smallabs{\sigma})^{-1-2\delta}\leq \rho^{-1-\delta+\alpha} \smallabs{\sigma}^{-\delta-\alpha}
$$
allows us to bound the RHS by
$$
|\sigma|^{-\alpha-\delta} \norm{\tP(\sigma) u}_{s,1/2-\alpha}.
$$
Meanwhile, the LHS is clearly larger than
$\norm{u}_{s,-3/2-\delta},$ and the result follows.
  \end{proof}

We will also require a slightly different special case of
Theorem~\ref{theorem:lowerbnd}, which we record for later use as a
separate lemma.
\begin{lemma}\label{lemma:bigf2}
  For all $\delta>0$ and $s>1+\delta,$
  $$
\tP(\sigma)^{-1}: H_b^{s,-1/2-\delta} \to |\sigma|^{-1-2\delta} H_b^{s,-3/2-\delta}
  $$
  \end{lemma}
  \begin{proof}
We take $l=-3/2-\delta$ in Theorem~\ref{theorem:lowerbnd}, which entails
$s>1/2+\delta$ and $\nu \in (-1-\delta, -\delta).$  Taking
$\nu=-2\delta,$ and estimating $(\rho+\smallabs{\sigma})^{-1-2\delta}
\leq \smallabs{\sigma}^{-1-2\delta}$ on the RHS gives the desired estimate.
   \end{proof}

Recall that
\begin{equation}\label{tPsplit}
\tP(\sigma) = \tP(0)-\sigma R
\end{equation}
where
\begin{equation}\label{R}
R= -2i\rho(\rho\partial_\rho-1)+\A^{\kappa+1}\Diffb^1(X) + \sigma \A^\kappa.
\end{equation}

In the following calculations we make use of the fact that $H_b^{s_1,\gamma_1} \subset H_b^{s_2,\gamma_2}$ for $s_2 \le s_1$ and $\gamma_2 \le \gamma_1$. From \eqref{R} we see
$$
R(\rho^n Y_\ell) = C_n \rho^{n+1} Y_\ell+ \A^{\kappa+n+1}+ \sigma
\A^{\kappa+n},
$$
where we crucially note that $C_1=0.$
Likewise,
\begin{equation}\label{Rmapping}
R(H_b^{s,\gamma}) \subset H_b^{s-1,\gamma+1} + \sigma H_b^{s,\gamma+\kappa}
\quad \text{ and } \quad
R(\A^\kappa) \subset \A^{\kappa+1} + \sigma \A^{2\kappa}.
\end{equation}
In consequence, the foregoing lemmas imply 
\begin{equation} \label{RPinv1} \begin{aligned}
&R\tP(0)^{-1}: H_b^{s,\gamma} \to &\sum_{j=3}^\infty
                 \rho^jY_{j-2} + H_b^{s+1,\gamma-1-} + \sigma
                                                         H_b^{s+2,\kappa-\frac{1}{2}-},\\
                                                         & &\gamma \in
                                                         \left(\frac{3}{2},
                                                         \kappa +
                                                         \frac{3}{2}\right],\quad
                                                             s>0,
                                                           \end{aligned}
                                                           \end{equation}
\begin{equation} \label{RPinv2}\begin{aligned}
 &R\tP(0)^{-1}: H_b^{s,\gamma} \to & H_b^{s+1,\gamma-1-}+
                  \sigma H_b^{s+2,\gamma + \kappa-2-},\\ & &\gamma \in
                  \left(\frac{1}{2},\frac{3}{2}\right), \quad s> 1,
                  \end{aligned}\end{equation}
\begin{equation}
  \label{RPinv3}
  \begin{aligned}
    &R\tP(0)^{-1}: \rho^n Y_\ell \to &\rho^{n-1} Y_\ell + \sum_{j=3}^\infty
\rho^{j} Y_{j-2}+\A^{\kappa+2-}+\sigma \A^{\kappa+2},\\ & &n \geq 4,\
\ell \geq n-2
\end{aligned}\end{equation}\begin{equation}
\label{RPinv4}R\tP(0)^{-1}: \rho^3 Y_\ell \to \sum_{j=3}^\infty
\rho^{j} Y_{j-2}+\A^{\kappa+2-}+\sigma \A^{\kappa+2},\  \ell \geq 1.
\end{equation}

Note that we have written sums of $\rho^jY_{j-2}$ terms with sums going out to
infinity for simplicity in bookkeeping, but all but a finite number of these
terms are subsumed in the conormal errors that we also carry
along. When $\gamma \in (\frac{1}{2},\frac{3}{2})$ we use Corollary
\ref{cor:prelimmapprop} to find $\tP(0)^{-1}: H_b^{s,\gamma} \to
H_b^{s,\gamma-2-}$ for $s > \frac{3}{2}-\gamma$. 

Additionally, Lemma~\ref{lemma:bigf2} and \eqref{R} yield
 the following estimate when we replace $\tP(0)$
with the full $\tP(\sigma)$ (and a factor of $\sigma$ thrown in for
purposes of later iteration):
\begin{equation}
\tP(\sigma)^{-1}(\sigma R):  H_b^{s,-3/2-\delta} \to
               |\sigma|^{-2\delta}H_b^{s-1,-3/2-\delta},\quad
                s>2+\delta.
\end{equation}
Since $[\sigma \pa_\sigma, \tP(\sigma)]$ is almost but not exactly $-\sigma R$, this is not quite the estimate we
will need
to obtain mapping properties of $(\sigma\pa_\sigma)^J
\tP(\sigma)^{-1}$ below; rather we need a slight variant to take into
account iterated commutators of $\sigma \pa_\sigma$ and $\tP(\sigma).$
To this end, recall that              $$
\tP(\sigma)= \tP(0)-\sigma R
$$
where $R=R_0+\sigma \A^\kappa,$
with $R_0$ independent of $\sigma.$ Thus
$$
[\sigma \pa_\sigma, \tP(\sigma)] = -\sigma R+\sigma^2 \A^\kappa,
$$
and moreover all iterated commutators of $\sigma \pa_\sigma$ with
$\tP(\sigma)$ are of the form
\begin{equation}\label{itcomm1}
[\sigma\pa_\sigma, -\sigma R+\sigma^2 \A^\kappa]=-\sigma R+\sigma^2 \A^\kappa.
\end{equation}
We thus remark for later use that changing $\sigma R$ by a
multiplication operator in $\sigma^2 \A^\kappa$ does not change its
mapping properties, since the mapping properties \eqref{Rmapping}
apply a fortiori when $R$ is replaced by a multiplication operator in
$\sigma \A^\kappa$ (such a term figured as part of $R$ itself in the
first place).  Thus, more generally,
\begin{equation}\label{perturbedR}
\tP(\sigma)^{-1}(\sigma R+\sigma^2\A^\kappa):  H_b^{s,-3/2-\delta} \to
                |\sigma|^{-2\delta}H_b^{s-1,-3/2-\delta},\quad
                s>2+\delta.
\end{equation}


\subsection{Conormal estimates} \label{subsection:conormalestimates}

              We begin by stating the low-frequency estimates on the
              twisted resolvent that are essential to our energy decay
              results. We will employ the notation
            $\fracpart{\kappa}\in [0, 1)$ for the fractional
part of $\kappa$ and $\floor{\kappa}\in \ZZ$ for the floor function, so that
$$
\kappa=\fracpart{\kappa}+\floor{\kappa}.
$$

              \begin{proposition}\label{prop:lowenergy}
Let $s>0.$
                The twisted resolvent enjoys the following
                low-frequency asymptotics: If $f \in  H_b^{s,\kappa+3/2}$, then
$$
 (\sigma\pa_\sigma)^M
\tP(\sigma)^{-1} f\in \smallabs{\sigma}^{1+\kappa-} L^\infty
 ((-1,1)_\sigma; H^2_{\loc}) +\CI((-1,1)_\sigma; H^2_{\loc})
$$
for
$$
M\leq s+\floor{\kappa}.
$$
                \end{proposition}

                \begin{proof}
                  Suppose $\tP(\sigma) u=f.$  
We approximate the solution by applying a formal Neumann series
argument: using the decomposition \eqref{tPsplit}, we have for all $N,$
\begin{equation}
  \label{eq:1}
  \begin{aligned}
    (\tP(0)-\sigma R)^{-1} &= \big((\Id-\sigma R \tP(0)^{-1}) \tP(0)\big)^{-1}\\
    &= \tP(0)^{-1}\big( \Id +\dots + (\sigma \tP(0)^{-1} R)^N \big)\\  
    	&\quad + \tP(0)^{-1}(\Id-\sigma  R\tP(0)^{-1})^{-1}(\sigma R
    \tP(0)^{-1} )^{N+1}\\
    &= \tP(0)^{-1}\big( \Id +\dots + (\sigma R \tP(0)^{-1} )^N \big)\\
    	&\quad +  (\tP(0)-\sigma  R)^{-1} (\sigma R \tP(0)^{-1} )^{N+1} .
  \end{aligned}
  \end{equation}
Define $f_n := (R\tP(0)^{-1})^n f$ so that
 	\[ u = \tP(0)^{-1}f_0 + \cdots + \sigma^N\tP(0)^{-1}f_N + \sigma^{N+1}\tP(\sigma)^{-1}f_{N+1}. \]
By \eqref{RPinv1},
	$$
		f_1 = R\tP(0)^{-1}f \in \sum_{j=3}^{\infty} \rho^{j}Y_{j-2} + H_b^{s+1,\kappa+\frac{1}{2}-} + \sigma H_b^{s+2,\kappa-\frac{1}{2}-}.
                $$
                (Recall that we will use $Y_j$ to denote finite linear
                combination of spherical harmonics of the given
                weight, without changing notation for each occurrence.)

Then \eqref{RPinv2}, \eqref{RPinv3}, \eqref{RPinv4} yield\footnote{We let the upper index of the
  sum of $\rho^j$ terms equal infinity for brevity; in fact of
  course we could rewrite this as a finite sum, with all terms beyond
  $j=\floor{\kappa}+2$ being absorbed in the $H_b^{s+2,\kappa - \frac{1}{2}-}$  term, with the caveat that polynomial dependence on $\sigma$ must then be allowed in that term.}
	\begin{align*}
        &f_2 = R\tP(0)^{-1}f_1 \\
        	&=\sum_3^\infty \rho^j(Y_{j-2} +\sigma Y_{j-2}+ Y_{j-1})
+ H_b^{s+2,\kappa-\frac{1}{2}-}+ \sigma H_b^{s+3,\kappa-\frac{3}{2}-} + \sigma^2 H_b^{s+4,\kappa-\frac{1}{2}-}.
\end{align*}
We will frequently be faced with terms that have polynomial dependence
on $\sigma,$ and will not be especially interested in the degrees of
the resulting polynomials (which could be bounded in terms of $\kappa$ in the iteration below if desired).  To streamline the
resulting notation, we therefore write
$$
\poly Y_j,\ \poly H_b^{s,l}
$$
to indicate respectively polynomials in $\sigma$ with coefficients in
$Y_j$ or $H_b^{s,l}.$  In this notation, then, 
we continue applying $R \tP(0)^{-1}$ to establish inductively that
	$$
		f_n \in \sum_{j=3}^\infty \rho^j
                \left(\sum_{\ell=j-2}^{j+n-3} \poly Y_\ell \right) +
                H_b^{s+n,\kappa+\frac{3}{2}-n-}+\sigma \poly H_b^{s+1+n,
                  \kappa+1/2-n-}
	$$
        for $1 \le n \le \floor{\kappa}$.
Hence, absorbing all terms in the first sum beyond $j=3$ in the following term,
$$        f_{\floor{\kappa}} \in \rho^3\left(
                  \sum_{\ell=1}^{\floor{\kappa}} \poly Y_\ell \right)
                + \mathbb{C}[\sigma] H_b^{s+\floor{\kappa},
                  \fracpart{\kappa}+\frac{3}{2}-}+
\sigma\poly H_b^{s+\floor{\kappa}+1, \fracpart{\kappa}+1/2-}.$$
                
We can continue one more step with the iteration: setting
        $$
J=\floor{\kappa}+1,
$$
we obtain
	$$
		f_{J} \in
                \mathbb{C}[\sigma] H_b^{s+\floor{\kappa}+1,\fracpart{\kappa}+\frac{1}{2}-}+
                \sigma \poly
                H_b^{s+\floor{\kappa}+2,\fracpart{\kappa}-\frac{1}{2}-}
                +\sigma^2 \poly H_b^{s+3+\floor{\kappa}, \fracpart{\kappa}+\frac{1}{2}-}
                  $$
where the terms containing the spherical harmonics $Y_\ell$ have been absorbed into the first term in $f_J$. We split this term into pieces
                  $$
f_{J}=F_{J}+\sigma G_J
$$
$$
F_J \in \mathbb{C}[\sigma] H_b^{s+\floor{\kappa}+1,\fracpart{\kappa}+\frac{1}{2}-}+\sigma^2 \poly H_b^{s+3+\floor{\kappa}, \fracpart{\kappa}+\frac{1}{2}-},
                  $$
                  $$
G_J \in \poly H_b^{s+\floor{\kappa}+2,\fracpart{\kappa}-\frac{1}{2}-}.
                  $$
                  Then
 \begin{equation}\label{mess}
                  \begin{aligned}
                    \tP(\sigma)^{-1} f &= \tP(0)^{-1} (f_0+\dots+\sigma^{J-1}
                    f_{J-1}) + \sigma^J \tP(\sigma)^{-1}
                    (F_J+\sigma G_J)\\
                   & =\tP(0)^{-1} (f_0+\dots+\sigma^{J-1}
                    f_{J-1}) + \sigma^J \tP(0)^{-1} F_J\\&\quad +\sigma^{J+1}
                    \tP(\sigma)^{-1}(R \tP(0)^{-1} F_J)  + \sigma^J
                    \tP(\sigma)^{-1} (\sigma G_J)\\
                    & =\tP(0)^{-1} (f_0+\dots+\sigma^{J-1}
                    f_{J-1}) + \sigma^J \tP(0)^{-1} F_J\\
                    &\quad +\sigma^{J+1}
                    \tP(\sigma)^{-1}(R \tP(0)^{-1} F_J  +G_J)
                                          \end{aligned}
                             \end{equation}
 where we have applied \eqref{eq:1} with $N=0$ in the penultimate
 step.  Note that the terms with polynomial dependence in $\sigma$
 will be favorable for obtaining regularity of the resolvent at
 $\sigma=0$: it is the final term
 $$
                    \tP(\sigma)^{-1}(R \tP(0)^{-1} F_J  +G_J)\equiv                     \tP(\sigma)^{-1}W
                    $$
                    with
$$W=                    (R \tP(0)^{-1} F_J  +G_J)$$
                    that will require finer analysis.

We now analyze the regularity of the terms above. For $n=1,\dots, J-1=\floor{\kappa},$ we have
	$$
\tP(0)^{-1} f_n =\sum_{j=1}^\infty \rho^j \sum_{\ell=j-1}^{j+n-1}
\poly Y_\ell + H_b^{s+n+2,\kappa-\frac{1}{2}-n-} + \sigma \poly H_b^{s+1+n,\kappa-n-\frac{3}{2}-},
	$$
using Lemmas \ref{lemma:P0inverse}, \ref{lemma:P0inverselessdecay},
and \ref{lemma:P0inversesh}.
Likewise, by Lemmas~\ref{lemma:P0inverse} and \ref{lemma:P0inverselessdecay},
$$
\tP(0)^{-1} F_J \in \poly H_b^{s+\floor{\kappa}+3, \fracpart{\kappa}-3/2-},
$$
while using \eqref{RPinv2} gives
	$$
		R\tP(0)^{-1}F_J  \in \poly H_b^{s+\floor{\kappa}+2, \fracpart{\kappa}-1/2-}.
	$$
In particular, we have now established
$$
W\in \poly H_b^{s+\floor{\kappa}+2, \fracpart{\kappa}-1/2-}.
$$
hence Lemma \ref{lemma:Psiginvlargeinput} yields
\begin{equation}\label{bigterms1}
  \tP(\sigma)^{-1}W  \in
\smallabs{\sigma}^{-1+\fracpart{\kappa}-}
L^\infty_\sigma H_b^{s+\floor{\kappa}+2, -\frac{3}{2}-}.
\end{equation}
where $  L^\infty_\sigma H_b^{\cdot,\cdot}$ denotes a bounded function
of $\sigma\in (-1,1)$ with values in the given Sobolev space.
We now sharpen this to obtain a partial conormality in $\sigma$ at
$\sigma=0:$ in particular, we claim
\begin{equation}\label{ourclaim}
(\sigma \partial_\sigma)^M \tP(\sigma)^{-1}W  \in 
 \smallabs{\sigma}^{-1+\fracpart{\kappa}-}
  L^\infty_\sigma H_b^{s+\floor{\kappa}+2-M, -\frac{3}{2}-}
\end{equation}
for all integers $$M <s+\floor{\kappa}+1.$$
To accomplish this, we will need a small result about commutators of
$\sigma \pa_\sigma$ with $\tP(\sigma)^{-1}.$  In what follows, we use
the letter $Q$ to denote an operator of the form
 $(\text{constant})\cdot (-\sigma R+\sigma^2 \A^k)$ but with the specific operator
allowed to change in each occurrence.
\begin{lemma}\label{lemma:commutatorwithinverse}
For all $M \in \NN,$ there exist constants $C_{\ell m}$ such that
  $$
(\sigma \pa_\sigma)^M \tP(\sigma)^{-1}=\sum_{\ell+m\leq M}C_{\ell m} (\tP(\sigma)^{-1} Q)^\ell
\tP(\sigma)^{-1} (\sigma\pa_\sigma)^m.
  $$
\end{lemma}
(Note that some of the $C_{\ell m}$ are in fact zero.)
\begin{proof} By induction on $M$, using the crucial fact \eqref{itcomm1}
that iterated commutators of $\sigma \pa_\sigma$ with $\tP(\sigma)$
all have the form of $Q.$\end{proof}

Now for any $m,$ Lemma \ref{lemma:Psiginvlargeinput} yields
$$
\tP(\sigma)^{-1} (\sigma \pa_\sigma)^m W \in \smallabs{\sigma}^{-1+\fracpart{\kappa}-}
L^\infty_\sigma H_b^{s+\floor{\kappa}+2, -\frac{3}{2}-}
$$
(since $\sigma \pa_\sigma$ passes harmlessly through a $\poly$
factor).  Hence repeated use of \eqref{perturbedR}
yields for $\ell \leq M,$
$$
(\tP(\sigma)^{-1} Q)^\ell
\tP(\sigma)^{-1} (\sigma\pa_\sigma)^m W  \in
 \smallabs{\sigma}^{-1+\fracpart{\kappa}-}
  L^\infty_\sigma H_b^{s+\floor{\kappa}+2-\ell, -\frac{3}{2}-};
  $$
the constraint on $M$ follows is the requirement that the
b-regularity index, which is $s+\floor{\kappa}+2-(M-1)$ after the
application of of $(\tP(\sigma)^{-1} Q)^{M-1}$ still remain greater
than $2$ as required for the $M$'th and final application of \eqref{perturbedR}.
Thus by Lemma~\ref{lemma:commutatorwithinverse} we obtain \eqref{ourclaim}.

Since the terms other than $\tP(\sigma)^{-1} W$ in \eqref{mess} are
smooth (indeed, polynomial) in $\sigma,$
we finally arrive at the estimate
\begin{equation}\label{twistedresolventbound}
f \in H_b^{s,\kappa+3/2-} \Longrightarrow (\sigma\pa_\sigma)^M
\tP(\sigma)^{-1} f\in \smallabs{\sigma}^{1+\kappa-} L^\infty
 ((-1,1)) +\CI((-1,1))
\end{equation}
(implicitly with values in the space $\Hb^{s+\floor{\kappa}+2-M, -3/2-}$)
for
\begin{equation}\label{M}
M\leq s+\floor{\kappa}.
\end{equation}
Since our estimates are with respect to the spaces
$\Hb^{s+\floor{\kappa}+2-M, -3/2-}$, the constraint \eqref{M} ensures
that these lie in $H^2_{\loc}.$
\end{proof}

The following lemma allows us to estimate $b$ regularity of data with
the oscillatory factor $e^{-i\sigma/\rho}$ inserted.
\begin{lemma}\label{lemma:conjugatedregularity}
Let $s \geq 0.$  If $f \in H_{b}^{s,\gamma+s}$ then $e^{-i\sigma/\rho} f \in H_b^{s,
  \gamma},$
uniformly for $\sigma$ in a compact set.
\end{lemma}
\begin{proof}
  We compute
$$\norm{(\rho\pa_\rho)^\alpha \pa_\theta^\beta e^{-i\sigma/\rho} f}=\norm{(\rho\pa_\rho-i\sigma/\rho)^\alpha \pa_\theta^\beta f}.$$
As long as $\alpha+\abs{\beta}\leq m,$ then, the differential operator
in question is in $\rho^{-m} \Diff_{b}^{m}$ (with coefficients
uniformly bounded for bounded $\sigma$) hence is the RHS is bounded by the
$H_{b}^{m,\gamma+m}$ norm. Thus we have obtained the result for integer $s.$  The
general case follows by interpolation.
\end{proof}
We may thus translate our estimates back to the setting of the
ordinary (unconjugated) resolvent:
\begin{corollary}\label{cor:resolventbound0}
Let $s > 0.$
  The resolvent enjoys the following low-frequency asymptotics: If $f \in  H_{b}^{s,\kappa+3/2+s}$ then
$$
 (\sigma\pa_\sigma)^M
P_\sigma^{-1} f\in \smallabs{\sigma}^{1+\kappa-} L^\infty
 ((-1,1)_\sigma; H^2_{\loc}) +\CI((-1,1)_\sigma; H^2_{\loc})
$$
for
$$
M\leq s+\floor{\kappa}.
$$
\end{corollary}
\begin{proof}
Note that
  $$
P_\sigma^{-1} f= e^{i\sigma/\rho} \tP(\sigma)^{-1} e^{-i\sigma/\rho}f.
$$
The leading factor $e^{i\sigma/\rho}$ is smooth in $\sigma \in (-1,1)$
uniformly on compact sets in $X^\circ,$ hence it suffices to verify
that $e^{-i\sigma/\rho}f$ satisfies the hypotheses of Proposition~\ref{prop:lowenergy}.
 This in turn follows from
Lemma~\ref{lemma:conjugatedregularity}.

\end{proof}

\section{Proof of Main Theorem} \label{MainProof}

   Recall by Proposition~\ref{prop:resex} that if $u$ solves the initial
value problem, then
	$$
		(2\pi)^{1/2}\check{u}(\sigma,\bullet) = P_\sigma^{-1}(-i\sigma u_0 + P_1 u_0 - u_1) \equiv P_\sigma^{-1}(\sigma f_0 +g_0) 
                $$
                with
                \begin{equation}\label{f0g0}
f_0=-i  u_0,\quad g_0 =P_1 u_0-u_1.
\end{equation}
We take the Fourier transform to recover $u$ and split the solution into low and high frequency parts (denoted $u_L$ and $u_H$, respectively). Let $\chi_{<1}(|\sigma|)$ be a cutoff function equal to 1 on $(-1/2, 1/2)$ and supported in $(-1,1)$ and take $\chi_{>1}(|\sigma|) = 1- \chi_{<1}(|\sigma|)$. We write $u = u_L + u_H$ where
\begin{equation}\label{uL}
	u_L = \int_{-\infty}^\infty \chi_{<1}(\sigma) e^{-i\sigma
          t}P_\sigma^{-1}(\sigma f_0 + g_0) \, d\sigma
\end{equation}
and
\begin{equation}\label{uH}
	u_H = \int_{-\infty}^\infty \chi_{>1}(\sigma) e^{-i\sigma t}P_\sigma^{-1}(\sigma f_0 + g_0) \, d\sigma.
\end{equation}

It now suffices to treat the asymptotic behavior of the high and low-frequency contributions separately.
We begin with $u_L.$

By Lemma~\ref{lemma:conormalFT} in Appendix~\ref{app:FT} and Sobolev embedding, it will suffice in estimating $u_L$ to show that
\begin{equation}\label{lowconormal}\check{u}_L(\sigma,x) \in \smallabs{\sigma}^{\kappa+1-} I^M L^\infty_{\loc}+\CI\end{equation} for some $M\geq
\kappa+1$
(where all spaces in $\sigma$ are valued in $H^2_{\loc}$).

To this end, we begin by noting (with a view to potential future
applications) the sharp hypotheses that are
necessary to obtain the estimate \eqref{lowconormal} for the
low-frequency part of the solution: what we will in fact use
is
\begin{equation}\label{sharplowhypoth}
  \begin{aligned}
    u_0 &\in \Hb^{3, \kappa+7/2},\\
    u_1 &\in \Hb^{2, \kappa+7/2}.
        \end{aligned}
\end{equation}

Now apply Corollary~\ref{cor:resolventbound0} with
$s=2$ to $$g_0=P_1 u_0 -u_1\in
\Hb^{2, \kappa+7/2} $$ to obtain
$$
(\sigma\pa_\sigma)^M
P_\sigma^{-1} g_0\in \smallabs{\sigma}^{1+\kappa-} L^\infty
 ((-1,1)_\sigma; H^2_{\loc}) +\CI((-1,1)_\sigma; H^2_{\loc}),\quad M
 \leq 2+\floor{\kappa}.
$$

Apply Corollary~\ref{cor:resolventbound0} with
with $\kappa$ replaced by $\kappa'=\kappa-1$ (hence the hypotheses on the
perturbation are a fortiori satisfied for $\kappa'$) and with $s=3$ to $$\sigma f_0= -i\sigma u_0\in
\sigma \Hb^{3, \kappa'+9/2} $$
to obtain for $M
\leq  3+\floor{\kappa'}$
$$
(\sigma\pa_\sigma)^M
P_\sigma^{-1} \sigma f_0 \in \smallabs{\sigma}^{2+\kappa'-} L^\infty
 ((-1,1)_\sigma; H^2_{\loc}) +\CI((-1,1)_\sigma; H^2_{\loc}),
 $$
 i.e. for $M \le 2 + \floor{\kappa}$ \,
 $$
(\sigma\pa_\sigma)^M
P_\sigma^{-1} \sigma f_0 \in \smallabs{\sigma}^{1+\kappa-} L^\infty
 ((-1,1)_\sigma; H^2_{\loc}) +\CI((-1,1)_\sigma; H^2_{\loc}).
 $$

By Lemma~\ref{lemma:conormalFT}, as noted above (and since
$2+\floor{\kappa} >1+\kappa$), this concludes the proof that $u_L$ has the desired decay, and we now
turn our attention to $u_H,$ the high-frequency component of the solution.

We will decompose the expression $P_\sigma^{-1}(\sigma f_0 + g_0)$ in
\eqref{uH} via an iterative argument. We approximate
$P_{\sigma}^{-1}(\sigma f_0 +g_0) \approx \sigma^{-1}f_0$ and let
$v_1$ denote the error. Recall $f_0 \in \Hb^{s+1,\kappa+\frac{7}{2}}$ and $g_0 \in \Hb^{s,\kappa+\frac{7}{2}}$. Direct calculation shows
		\begin{align*}
			P_\sigma v_1 =  (g_0 - iP^1f_0) + \sigma^{-1}(\Delta + P^2)(-f_0)  =: f_1 + \sigma^{-1}g_1
		\end{align*}
	where $f_1 \in H_b^{s,\kappa+\frac{7}{2}}$ and $g_1 \in H_b^{s-1,\kappa+\frac{11}{2}}$ (since $P^1: \Hb^{m,\ell} \to \Hb^{m-1,\ell+\kappa+1}$ and $(\Delta+P^2): \Hb^{m,\ell} \to \Hb^{m-2,\ell+2}$). Now we have
	\[ P_\sigma^{-1}(\sigma f_0 + g_0) = \sigma^{-1}f_0 + P_\sigma^{-1}(f_1 + \sigma^{-1}g_1 ). \]
Next we iterate the process and approximate $P_\sigma^{-1}(f_1 + \sigma^{-1}g_1 ) \approx \sigma^{-2}f_1$: 
	\begin{align*} 
		P_\sigma^{-1}(f_1 + \sigma^{-1}g_1 ) &= \sigma^{-2}f_1 + P_\sigma^{-1}(\sigma^{-1}(g_1-iP^1f_1) +\sigma^{-2}(\Delta+P^2)(-f_1) )\\
			&=: \sigma^{-2}f_1 + P_\sigma^{-1}(\sigma^{-1}f_2 +\sigma^{-2}g_2 )
	\end{align*}
where $f_2 \in H_b^{s-1,\kappa+\frac{11}{2}}$ and $g_2 \in H_b^{s-2,\kappa+\frac{11}{2}}$. Thus we have

	\begin{equation} \label{rsigmaiteration}
		P_\sigma^{-1} (\sigma f_0 + g_0) = \underbrace{ \sigma^{-1} f_0 + \sigma^{-2}f_1}_{=:\check{u}_a} + \underbrace{P_\sigma^{-1}\left(  \sigma^{-1} f_2 + \sigma^{-2} g_2 \right)}_{=: \check{u}_b} 
	\end{equation}


We plug \eqref{rsigmaiteration} into our expression for $u_H$ (see equation \eqref{uH}) and bound
each term separately. Note $|f_0|+ |f_1| \lesssim \langle r \rangle^{-5-\kappa}$ so we calculate for any $M \ge 1$:
	\begin{align*}
		&\left|\int_{\sigma \in \R} \chi_{>1}(|\sigma|) \check{u}_a(\sigma)e^{-it\sigma} \dd \sigma \right| \\		
		&\lesssim  \la r \ra^{-\kappa -5}\la t \ra^{-M} \left| \int \partial_{\sigma}^M\big[ \chi_{>1}(|\sigma|)( \sigma ^{-1} + \sigma^{-2})\big] e^{-i\sigma t} \dd \sigma \right| \\
		&\lesssim \la t \ra^{-M} \la r \ra^{-\kappa-5}.
	\end{align*}
	
	By Proposition \ref{largesigmaresbnd}, we have $ |(\sigma \partial_{\sigma})^{M} (\check{u}_b(\sigma)e^{-i\sigma  r })| \lesssim |\sigma|^{M-2}\la r \ra^{-1+M(1-\varepsilon}$ for $s - 2 \ge (\loss +1)(2M + 5)$ and $M < \kappa+5$ with $\varepsilon \in (0,1]$. We use this to calculate
	\begin{align*}
			&\left|\int_{\sigma \in \R} \chi_{>1}(\sigma) \check{u}_b(\sigma)e^{-it\sigma} \dd \sigma \right|\\ 				
				&\approx \la t-r \ra^{-M} \left| \int \sigma^{-M} \left[ \sum_{\ell=0}^M (\sigma\partial_{\sigma})^{\ell}  (\chi_{>1}(|\sigma|)\check{u}_b e^{-ir\sigma}) \right] e^{i(r-t)\sigma} \dd \sigma \right| \\
				&\lesssim \la r \ra^{-1+M(1-\varepsilon)} \la t -r \ra^{-M}.
		\end{align*}

	Combining the above results, we find 
	\[ |u_H(t,x)| \lesssim \la t \ra^{-M} \la r \ra^{-\kappa-1} + \la r \ra^{-1+M(1-\varepsilon)} \la t -r \ra^{-M}. \]
The main theorem (and indeed a finer estimate near the light cone for
this part of the solution) then follows in the high frequency case when we take $M= \kappa + 2 - \epsilon$ (and thus the b-regularity requirement is $s > (\loss+1)(2\kappa + 9) + 2$).\qed

    \appendix

\section{Mellin transforms and b-Sobolev spaces}\label{app:bgeometry}
In this appendix, we describe the Mellin transform characterization of
the b-Sobolev spaces defined in Section~\ref{section:background} above.

   If $u$ is a function of $\rho, \theta$ ($\rho \in [0, \infty)$)
set
$$
\M u (\xi,\theta) = \int_0^\infty \rho^{-i\xi} u(\rho,\theta) \frac{d\rho}\rho.
$$

Let $\hol(\Omega)$ denote the space of holomorphic functions on the domain $\Omega \subset \CC$ with values in $L^2(S^2).$

Let $L^2_{bb}$ denote the space of $u$ that are $L^2$ with respect to the ``b''-density $d\rho/\rho$.  Let $L^2_{bb,A}$ denote the subspace of functions supported in $\rho \in [0,A].$
\begin{lemma}
  The range of $\M$ on $L^2_{bb, A}$ is
  $$
\big\{ w(\xi,\theta)\colon w \text{ holomorphic in } \Im \xi>0,\ \sup_{\mu>0} A^{2\mu}\int_{\Im \xi=\mu} \smallnorm{w(\xi,\theta)}_{L^2(S^2)}^2\, d\xi<\infty\big\}.
  $$
    \end{lemma}
    \begin{proof}
      Set $x=\log \rho.$  Then
      $$
\M u(\xi) = \int_{-\infty}^\infty  e^{-i\xi x} u(e^x,\theta) \, dx=(2\pi)^{1/2} \F (u \circ \exp) (\xi),
$$
with $u \circ \exp\in L^2(\RR_x \times S^2_\theta)$ and supported in $x\in (-\infty, \log A].$
Let us take $A=1$ for now, hence the support of $u \circ \exp$ is $(-\infty, 0].$
The Paley--Wiener theorem tells us that the Fourier transforms of such functions are precisely the space
$$
\big\{ w(\xi,\theta)\colon w\in \hol( \Im \xi>0),\ \sup_{\mu>0} \int_{\Im \xi=\mu} \smallnorm{w(\xi,\theta)}_{L^2(S^2)}^2\, d\xi<\infty\big\}.
$$
This establishes the result for $A=1.$  More generally,  $u \in L^2_{bb,A}$ iff $u(A^{-1} \rho, \theta)\in L^2_{bb,1},$ so the result follows from the fact that $\M [u(A^{-1} \rho, \theta)](\xi)=A^{i\xi} \M u.$
      \end{proof}

      Now we observe that adding a weight simply shifts the domain of holomorphy, since $\M \rho^\alpha u(\xi)=\M u (\xi+i\alpha)$ hence
      \begin{align*}
&\M \rho^{\alpha} L^2_{bb, A}\\
&\to \big\{ w(\xi,\theta)\colon w \in \hol(\Im \xi>-\alpha),\ \sup_{\mu>-\alpha} A^{2\mu}\int_{\Im \xi=\mu} \smallnorm{w(\xi,\theta)}_{L^2(S^2)}^2\, d\xi<\infty\}.
\end{align*}

Finally, we aim to keep track of $b$-regularity.  Recall that $H_b^{m}$ denotes
the (unweighted) $b$-Sobolev space of order $m,$ \emph{measured with respect to the
  metric volume form}, which is proportional to $d\rho/\rho^4\, d\theta;$ note
in particular that this entails the numerology 
$$
H_b^{0} = \rho^{3/2} L^2_{bb}.
$$
More generally, recall that $H_b^{m,\ell} = \rho^\ell H_b^{m,0}.$  We write $H^\bullet_{b,A}$ as before to denote the functions supported in $\rho\leq A.$

\begin{proposition} \label{prop:mellin}
For $m\geq 0$ the Mellin transform is an isomorphism
\begin{multline}\M: H_{b,A}^{m,\ell}\to \big\{ w(\xi,\theta)\colon w\in \hol( \Im \xi>-\ell-3/2),\\ \sup_{\mu>-\ell-3/2} A^{2\mu}\int_{\Im \xi=\mu} \smallnorm{(\xi^2+\Lap_\theta)^{m/2} w(\xi,\theta)}_{L^2(S^2)}^2\, d\xi<\infty\big\}.
\end{multline}  \end{proposition}
\begin{proof}
For $m$ an even integer, the result follows by the b-elliptic regularity of the differential operator $((rD_r)^2 +\Lap_\theta)^{m/2}.$  For general $m,$ it follows by interpolation and duality.
  \end{proof}

\section{Fourier transforms of finite-regularity conormal
  distributions}\label{app:FT}

We consider distributions in $\smallabs{\sigma}^\alpha I^mL_c^\infty(\RR; 0);$ this
is defined as the space of distributions on $\RR$ lying in
$\smallabs{\sigma}^\alpha L_c^\infty$ (with $c$ denoting compact support) and enjoying
$m$-fold iterated regularity under
vector fields tangent to the origin, which is to say, under powers up
to the $m$th of $\sigma
\pa_\sigma.$
We will
only deal with $m \in \NN$ in order to keep the discussion simple.
Thus $u\in \smallabs{\sigma}^\alpha I^m L_c^\infty(\RR;0)$ if $u$ is compactly
supported and enjoys the estimate
$$
\abs{\pa_\sigma^j u}\leq C_j \abs{\sigma}^{\alpha-j},\quad j=0,\dots, m.
$$
\begin{lemma}\label{lemma:conormalFT}
  Let $u \in \smallabs{\sigma}^\alpha I^m L_c^\infty(\RR;0),$ and assume
  $\alpha>0$ and $m\geq \alpha+1.$  Then
  $$
\F(u) =O(\ang{t}^{-1-\alpha}).
  $$
\end{lemma}
\begin{proof} (cf.\ Lemma 3.6 of \cite{Hi:20}.)
  We write
  $$
\F u=(2\pi)^{-1/2} \int u(\sigma) e^{i \sigma t}\, dt=w_0+w_1
$$
where
$$
w_0=(2\pi)^{-1/2} \int_{\smallabs{\sigma}<t^{-1}} u(\sigma) e^{i \sigma
  t}\, dt,\quad w_1=\F u-w_0.
$$
Then
$$
\abs{w_0}\lesssim \smallabs{t}^{-1-\alpha}
$$
since $u \in \smallabs{\sigma}^\alpha L^\infty.$  On the other hand,
integration by parts using the operator $(t^{-1}D_\sigma)^m$ in the
integral expression for $w_1$ yields a bulk term bounded by
$$
t^{-m} \int_{t^{-1}<\smallabs{\sigma}<C} \sigma^{\alpha-m}\, d\sigma
$$
(where $C$ depends on the support of $u$),
as well as boundary terms bounded by
$$
t^{-m} \sigma^{-\alpha-m'}\big\rvert_{\sigma=t^{-1}},\quad m'
=0, \dots, m-1.
$$
 Since $m\geq \alpha+1,$ these terms are all bounded by multiples of
 $\smallabs{t}^{-1-\alpha}$ for $\smallabs{t}$ large.
  \end{proof}

\bibliographystyle{alpha} 
        
        \bibliography{price}


\end{document}